\documentclass[11pt]{article}
\usepackage{graphicx}
\usepackage{amsmath}
\usepackage{amsfonts}
\usepackage{epsfig}
\usepackage{setspace}

\newcommand{\be}{\begin{equation}}
\newcommand{\ee}{\end{equation}}
\newcommand{\beqn}{\begin{eqnarray}}
\newcommand{\eeqn}{\end{eqnarray}}
\newcommand{\beqns}{\begin{eqnarray*}}
\newcommand{\eeqns}{\end{eqnarray*}}

\newcommand{\card}{\mbox{card}}
\newcommand{\Var}{\mbox{Var}}

\newcommand{\EE}{\ensuremath{{\mathbb E}}}
\newcommand{\II}{\ensuremath{{\mathbb I}}}

\newcommand{\fr}[1]{(\ref{#1})}

\newcommand{\Te}{\Theta}

\newtheorem{lemma}{Lemma}
\newtheorem{theorem}{Theorem}

\newtheorem{remark}{Remark}

\begin{document}

\title{\Large{\bf Minimax estimation for simultaneous deconvolution model in random design using wavelets}}

\author{
\large{Rida Benhaddou}  \footnote{E-mail address: Benhaddo@ohio.edu}
  \\ \\
Department of Mathematics, Ohio University, Athens, OH 45701} 
\date{}

\doublespacing
\maketitle
\begin{abstract}
Simultaneous deconvolution model is investigated when the design points are random and follow bounded densities and errors are either $i.i.d.$ or long-memory dependent. We propose wavelet hard-thresholding estimators, and show that they achieve asymptotically near-optimal convergence rates for the $L^2$-risk when the response function belongs to Besov space under both $i.i.d.$ errors and homoskedastic, long-memory errors. However, under heteroskedastic, long-memory errors the estimator is not optimal. For the $i.i.d.$ errors, the convergence rates are {identical} to those found in simultaneous deconvolution under regular design {and $i.i.d.$ errors}. However, for long-memory dependent errors, the convergence rates are new {and completely different from those in simultaneous deconvolution under regular design and long-memory errors}.   \\


{\bf AMS (2000) Subject Classification: 62G05, 62G20, 62G08 }
 \end{abstract} 

\section{Introduction}

Consider the problem of estimating a periodic function $f(t)$ based on observations  
\be
y_{li}=q_l(t_{li})+ \sigma(t_{li}) \varepsilon^{(l)}_{i}=\int^1_0f(s)g_l(t_{li}-s)ds+\sigma_l(t_{li})\varepsilon^{(l)}_{i}, \ \ \ l=1, 2, \cdots, {\bf r}, \ i=1, 2, \cdots, N, \label{conveq}
\ee
where, ${\bf r}$ is finite, and $t_{li}$ are $i.i.d.$ random variables with compactly supported probability density functions $h_l$. In addition, the design points $t_{li}$ are independent across the ${\bf r}$ channels. $\sigma_l(t)$ are known, positive variance functions. The function $f(t)$ is $L^2\left([0,1]\right)$, and the convolution kernels $g_l(t_l)$ are known periodic functions. Representation \fr{conveq} is called simultaneous deconvolution model or multichannel deconvolution model. Furthermore, $\{\varepsilon^{(l)}_{i}\}_{i\geq1}$ {are independent of $t_{li}$ and} are independent across the channels. Also, $\{\varepsilon^{(l)}_{i}\}_{i\geq1}$ are {(1) $i.i.d.$ zero-mean Gaussian sequences, or are (2)} zero-mean, stationary Gaussian sequences, for any  $i\in \{1, 2, \cdots, N\}$, exhibiting the structure
\be
\varepsilon^{(l)}_{i}=\sum^{\infty}_{m=0}a^{(l)}_m\eta^{(l)}_{i-m},\ with\ a^{(l)}_o=1,
\ee
where, $\{\eta^{(l)}_i\}_{i \in Z}$ are $i.i.d.$ Gaussian sequences that satisfy the property
\beqns
\lim_{m\rightarrow \infty} a^{(l)}_m m^{(\alpha_l+1)/2}=1,\ for\ \alpha_l \in (0, 1),
\eeqns
and
\be \label{longm1}
\Var\left(\sum^N_{i=1}\varepsilon^{(l)}_{i}\right)\sim \pi_{l} N^{2-\alpha_l},
\ee
where $\pi_{l}$ are finite positive constants, and $\alpha_l$ are long-memory parameters.  

Deconvolution is a common problem in image and signal processing. This includes, for instance, light detection and ranging (LiDAR) for remote sensing and reconstruction of blurred images. LiDAR is a method for determining ranges by targeting an object or a surface with a laser and measuring the time for the reflected light to return to the receiver. However, if the response function of the LiDAR is longer than the time resolution interval, the LiDAR signal may become blurry and the accuracy of the LiDAR deteriorates (see e.g. De Canditiis $\&$ Pensky~(2006)). This imprecision can be corrected using deconvolution techniques. In addition,  LiDAR may operate in a fixed direction or it may scan multiple directions. If ${\bf r}\geq 2$, when multiple LiDAR devices are used to recover a signal, the process is called a multichannel deconvolution or simultaneous deconvolution. 

Wavelet deconvolution has been studied extensively in the literature since the {seminal} work of Donoho~(1995a). The list includes Donoho~(1995b), Donoho and Johnstone~(1998), Abramovich and Silverman~(1998), Walter and Shen~(1999), Donoho and Raimondo~(2004), Johnstone et al.~(2004), to name a few. In the case of simultaneous deconvolution, the problem was investigated in De Canditiis and Pensky~(2006), while in functional deconvolution, it was studied in Pensky and Sapatinas~(2009, 2010, 2011), Benhaddou, Pensky and Picard~(2013), Benhaddou~(2017), Wishart~(2019), and Benhaddou and Liu~(2019).

Long memory has been investigated quite considerably in deconvolution models under regular design. The list of articles includes, for standard deconvolution, Wang~(1997), Wishart~(2013), and Benhaddou~(2016), for multichannel deconvolution, Benhaddou et al.~(2014) {where the number of channels increases as $N$ increases},  and Kulik, Sapatinas and Wishart~(2015) {where the number of channels is fixed}, and for functional deconvolution, Benhaddou and Liu~(2019). 

All the procedures mentioned above assume the sample design is regular, evenly spaced, but in reality this is usually not true. One of earliest to study wavelet estimation under irregular design  was Cai and Brown~(1998, 1999). This was followed by several articles, including in chronological order, Chesneau~(2007a, 2007b), Kulik and Raimondo~(2009), Pensky~(2013), Antoniadis et al.~(2014), Benhaddou~(2023, 2025), and Benhaddou and Liu~ (2024). 

{ In LiDAR applications, raw data points are inherently irregular. This can be seen in applications such as atmospheric research and environmental studies in which airborne LiDAR sensors are used. These variations in point density and spacing can be due to several factors, and the interested reader may refer to Petras et al.~(2023)}. 

The objective of the paper is to investigate simultaneous deconvolution model under {both $i.i.d.$ and} long-memory Gaussian errors via wavelet hard-thresholding procedures when the design points are random and follow bounded probability density functions $h_l$,  $l=1, 2, \cdots, {\bf r}$. To the best of our knowledge, this scenario has never been treated before. We would like to see how having long-memory errors, {both homoskedastic and heteroskedastic}, and random design points affect the convergence rates of the wavelet procedure under simultaneous convolutions setting. Under both $i.i.d.$ and homoskedastic long-memory errors, we construct an adaptive estimator { with respect to regularity of $f(t)$ and long-memory parameters} $\alpha_l$, $l=1, 2, \cdots, {\bf r}$, and {show} that it is asymptotically near-optimal {in the minimax sense} under  $L^2$-risk. However, under heteroskedastic long-memory errors, the estimator is not adaptive with respect to the long-memory parameters, nor is it optimal. In addition, under $i.i.d.$ errors, we notice the usual elbow effect that is often encountered in deconvolution problems under regular design and $i.i.d.$ errors. Furthermore, under homoskedastic long-memory errors, we see the usual elbow effects observed under the $i.i.d. $ errors in addition to a long-memory phase. {This latter} depends only on the strongest long memory across the channels. {This behavior is completely new, and it points to fact that randomization may not improve the performance of the estimation under long-memory errors in the case of simultaneous deconvolution when long memory is sufficiently strong}. Finally, under heteroskedastic long-memory errors, the thresholding estimator involves applying special weights to  the empirical wavelet and scaling coefficients {that are different from those found in Kulik, Sapatinas and Wishart~(2015) in their treatment of simultaneous deconvolution under regular, evenly spaced design and long-memory errors}. {Nevertheless}, our {upper bounds for the $L^2$-risk} {are slightly different from theirs when long-memory is strong across the channels}. 
\section{Wavelet Representation of Functions}
Consider a band-limited periodized wavelet basis $\psi$ (e.g. Meyer-type) on $[0, 1]$.  Let $m_{o}$ be the lowest resolution level and denote the scaling function by $\psi_{m_{o}-1, k}$. Define the set
\be
\Omega=\left\{\omega= (j, k):j=m_{o}-1, \cdots, \infty; k=0, 1, \cdots, 2^{j}-1\right\}. 
\ee
Since the functions $\{ \psi_{m_{o}-1, k}, \psi_{j, k}\}$ form an orthonormal basis of the $L^2([0, 1])$-space, the function $f(t)\in L^2([0, 1])$ admits the representation 
\be\label{invwavtrans}
f(t)=\sum_{\omega\in \Omega}\beta_{\omega} \psi_{\omega}(t),
\ee
where 
\be \label{wavtrans}
\beta_{\omega}=\int^1_0f(t)\psi_{\omega}(t)dt.
\ee
Denote the Fourier transform of $\psi_{\omega}(t)$ by $\tilde{\psi}_{\omega}(m)=\tilde{\psi}_{jk}(m)$. It is well known (see, e.g. Johnstone et al.~(2004), Section 3.1) that under the Fourier domain and for any $j>m_{o}$, one has 
\be \label{w-set}
W_{j}=\left\{m:\tilde{\psi}_{\omega}(m)\neq 0\right\}\subset \frac{2\pi}{3}\left[-2^{j+2}, -2^{j}\right]\cup \left[2^{j}, 2^{j+2}\right]. 
\ee
The quantity \fr{wavtrans} is referred to as wavelet transform (WT), and formula \fr{invwavtrans} is called the inverse wavelet transform (IWT). Under the regular design, observations are taken over a regular grid $i/N$, $i=1, 2, \cdots, N$, and both \fr{invwavtrans} and \fr{wavtrans} can be computed using standard algorithms. However, under random design, the computation  of these quantities requires a special treatment. \\
To gain some insight on how to handle randomness, let us consider the idealization of representation \fr{conveq}. Define the ${\bf r}$ random processes $Y_l$
\be \label{conveq-ideal}
Y_l(dt)=f*g_l(t)dt+\sigma_l(t)\varepsilon^{\alpha_l}W^{\alpha_l}_l(dt), \ \ l=1, 2, \cdots, {\bf r},
\ee
where $W^{\alpha_l}_l$ are independent fractional Wiener processes, and $\varepsilon^{\alpha_l}=N^{\alpha_l/2}$ (see, e.g. Benhaddou and Liu~\cite{ben-liu3}). The problem of solving the system \fr{conveq-ideal} for $f$ is well-posed as long as the Fourier transforms $\tilde{g}_l(m)$ of the functions $g_l(t)$ are entire functions on the complex $m$-plane, and 
\beqns
\sqrt{\sum^{\bf r}_{l=1}|\tilde{g_l}(m)|^2}\geq a e^{-2\pi b|\Im(m)|}(1+|m|)^{-L},
\eeqns
for some $a$, $b>0$, and $L\geq 1$ (see, e.g. De Canditiis and Pensky~(2006)). Then, there exist entire functions $\tilde{\mu}_l(m)$, $l=1, 2, \cdots, {\bf r}$, such that 
\beqns
\sum^{\bf r}_{l=1}\tilde{g_l}(m)\overline{\tilde{\mu}_l(m)}=1,
\eeqns
which implies 
\beqns
\sum^{\bf r}_{l=1}(g_l*\mu_l)(t)=\delta(t),
\eeqns
where $\delta(t)$ is the delta function, and $f$ can be recovered as 
\be
f(t)=(f*\delta)(t)=f*\sum^{\bf r}_{l=1}(g_l*\mu_l)(t)=\sum^{\bf r}_{l=1}(q_l*\mu_l)(t).
\ee
Thus, define the functions $U^{(l)}_{\omega}(t)$ such that
\be\label{u-f}
U^{(l)}_{\omega}(t)=\sum_{m\in W_{j}}\frac{\theta_l\tilde{g_l}(m)\overline{\tilde{\psi}_{\omega}(m)}}{\sum^{\bf r}_{l=1}\theta_l|\tilde{g_l}(m)|^2}e^{i2\pi mt}, 
\ee
where $\theta_l$ are weights that will be determined later. 
Notice that by Plancherel's identity and some properties of Fourier transform, the functions $U^{(l)}_{\omega}(t)$ defined in \fr{u-f}, satisfy the property
\beqn
\sum^{\bf r}_{l=1}\int^1_0{U^{(l)}_{\omega}(t)}q_l(t)dt &=&\sum^{\bf r}_{l=1}\int^1_0\sum_{m\in W_{j}}\frac{\theta_l\tilde{g_l}(m)\overline{\tilde{\psi}_{\omega}(m)}}{\sum^{\bf r}_{l=1}\theta_l|\tilde{g_l}(m)|^2}e^{-i2\pi mt}(f*g_l)(t)dt  \nonumber\\
&=&\sum^{\bf r}_{l=1}\sum_{m\in W_{j}}\frac{\theta_l\tilde{g_l}(m)\overline{\tilde{\psi}_{\omega}(m)}\tilde{f}(m)\overline{\tilde{g}_l(m)}}{\sum^{\bf r}_{l=1}\theta_l|\tilde{g_l}(m)|^2}  \nonumber\\
&=&\sum_{m\in W_{j}}{\overline{\tilde{\psi}_{\omega}(m)}}\tilde{f}(m)= \beta_{\omega}.
\eeqn
 Next is the list of conditions that will be used in the derivation of the theoretical results. \\
{\bf Assumption A.1.}\label{A1} $f \in {\bf{L^{2}}}\left[0, { 1} \right]$ is {uniformly bounded}, that is, there exists positive constant $M<\infty$ such that  $|f(t)|\leq M$, for all  $t\in \left[0, { 1} \right]$.\\
{\bf Assumption A.2.} \label{A2} The probability density functions $h_l$ are uniformly bounded {and bounded away from zero}, that is, on  $\left[0, {1} \right)$ there exist constants $m_{l1}$ and $m_{l2}$, with $0< m_{l1}\leq m_{l2} < \infty$, such that $m_{l1} \leq h_l(t) \leq m_{l2}$. 
 \begin{remark}	  
{ In LiDAR surveys applications, it is often encountered that the design density behaves like a {truncated normal}, or a {truncated exponential} distribution, which satisfy Assumption {\bf A.2}. }
  \end{remark}
{\bf Assumption A.3.}\label{A3}  The functions $\sigma_l(t)$ are uniformly bounded on $[0, 1]$, that is, there exist constants $0<\sigma_{lo}\leq\sigma_{l1}<\infty$ such that $\sigma_{lo}\leq \sigma_l(t)\leq \sigma_{l1}$.\\
{\bf Assumption A.4.}  \label{A4} For some positive constants $\nu_l$, $0<K_{l1}<K_{l2}<\infty$, independent of $m$, the Fourier coefficients $\tilde{g_l}(m)$ of the kernels $g_l(t)$ satisfy 
\be \label{gm}
K_{l1}|m|^{-2\nu_l}\leq \left|\tilde{g_l}(m)\right|^2\leq K_{l2}|m|^{-2\nu_l}. 
\ee
The parameters $\nu_l$ represent the degrees of ill-posedness of the deconvolution problem. \\
Define the quantity,
\be\label{Dell}
\Delta=\sum^{\bf r}_{l=1}\int^1_0{U^{(l)}_{\omega}(t)}\sigma_l(t_{l})dt.
\ee
 \begin{remark}	  
Notice that $\Delta$ defined in \fr{Dell}, will be zero whenever $\sigma_l(t_{l})$, $l=1, 2, \cdots, {\bf r}$, are constant, since Fourier basis functions are orthogonal to constants for all nonzero {frequencies}. {Therefore, under homoskedastic noise, $\Delta$ vanishes}.  $\Delta$ {would} not be zero when $\sigma_l(t_{l})$, $l=1, 2, \cdots, {\bf r}$, is non-constant, since
\beqn\label{delsigma}
\sum_{m\in W_j}\tilde{\psi}_{\omega}(m)\tilde{g_l}(m)\overline{\tilde{\sigma}_l(m)}\neq 0.
\eeqn
{Therefore, under heteroskedastic noise, $\Delta$ does not vanish}.
 \end{remark}
\section{Estimation algorithm}
{In {many} multi-LiDAR surveys, identical spatial distributions of points are unlikely to occur. This is due to several factors, such as environmental conditions, fog, dust, sunlight, etc, geometric discrepancies. {As a result}, achieving perfect alignment across all sensors is difficult. {Thus}, it is reasonable to consider $h_d\neq h_j$, with $d, j=1, 2, \cdots, {\bf r}$, $j\neq d$.} Indeed, consider the unbiased estimator 
\be \label{betest1}
\widehat{\beta_{\omega}}=\frac{1}{N}\sum^N_{i=1}\sum^{\bf r}_{l=1}\frac{y_{li}U^{(l)}_{\omega}(t_{li})}{h_l(t_{li})}.
\ee
Define the set 
\be\label{omega-j}
\Omega(J)=\left\{\omega= (j, k):m_{o}-1\leq j\leq J-1; k=0, 1, \cdots, 2^{j}-1\right\}. 
\ee
A hard thresholding estimator for $f$ is given by 
\be\label{fest}
\widehat{f}(t)=\sum_{\omega\in \Omega(J)}\widehat{\beta_{\omega}}\II\left(|\widehat{\beta_{\omega}}|> \lambda(j)\right)\psi_{\omega}(t),
\ee
where, the quantities $J$ and $\lambda(j)$ will be determined later. 
 \begin{lemma} \label{lem:Var-nonunif}
Let $\widehat{\beta_{\omega}}$ be defined as in \fr{betest1} and let assumptions ${\bf{A.1}}$-${\bf{A.4}}$ hold. {Assume $h_l(t)$ are known}, with $l=1, 2, \cdots, {\bf r}$. Under $i.i.d.$ errors or {long-memory errors} with $\Delta=0$, choose $\theta_l$ in \fr{u-f} as
\be\label{theta1-iid}
\theta_l=m_{l1}. 
\ee
Then, as $N \rightarrow \infty$, one has  
\beqn 
\EE\|\widehat{\beta_{\omega}}-\beta_{\omega}\|^2 &\asymp& \frac{\left[\sum^{\bf r}_{l=1}m_{l1}2^{-2\nu_l j}\right]^{-1}}{N}.  \label{var-bet-nonunif-iid}\\
\EE\|\widehat{\beta_{\omega}}-\beta_{\omega}\|^4 &\asymp&\frac{\left[\sum^{\bf r}_{l=1}m_{l1}2^{-2\nu_l j}\right]^{-2}}{N^2}.\label{var2-bet-nonunif-iid}
\eeqn
In addition, under {long-memory errors} with $\Delta\neq 0$, choose $\theta_l$ in \fr{u-f} as 
\be\label{theta1-lm-unif}
\theta_l=\left[\frac{1}{Nm_{l1}}+\frac{2^{-j}\pi_l}{N^{\alpha_l}}\right]^{-1}. 
\ee
Then, as $N \rightarrow \infty$, one has 
 \beqn 
\EE\|\widehat{\beta_{\omega}}-\beta_{\omega}\|^2 &\asymp&\left[\sum^{\bf r}_{l=1}\left[\frac{1}{Nm_{l1}}+\frac{2^{-j}\pi_l}{N^{\alpha_l}}\right]^{-1}2^{-2\nu_l j}\right]^{-1} .  \label{var-bet-nonunif-lm}\\
\EE\|\widehat{\beta_{\omega}}-\beta_{\omega}\|^4 &\asymp&  \left[\sum^{\bf r}_{l=1}\left[\frac{1}{Nm_{l1}}+\frac{2^{-j}\pi_l}{N^{\alpha_l}}\right]^{-1}2^{-2\nu_l j}\right]^{-2}.\label{var2-bet-nonunif-lm}
\eeqn
\end{lemma}
Let $\nu_o=\min\{\nu_1, \nu_2, \cdots, \nu_{\bf r}\}$, $\alpha_o=\min\{\alpha_1, \alpha_2, \cdots, \alpha_{\bf r}\}$, $\alpha^*=\alpha_{l^*}$, $\nu^*=\nu_{l^*}$, where
\be\label{l-optimal*}
  l^*=\arg \min_{1\leq l\leq {\bf r}}\left\{\frac{2^{2\nu_l}}{N^{\alpha_l}}\right\}.
\ee
According to {\bf Lemma} \ref{lem:Var-nonunif}, under $i.i.d.$ errors or {homoskedastic long-memory errors} ($\Delta=0$), the theoretical level-dependent thresholds in \fr{fest} should be chosen as 
\be\label{thresh-lm1-nonunif-iid}
\lambda(j)=\rho \frac{\left[\sum^{\bf r}_{l=1}2^{-2\nu_l j}\right]^{-0.5}{\ln^{\gamma}(N)}}{\sqrt{N}},
\ee
and, under {heteroskedasitic} long-memory errors ($\Delta\neq 0$), the thresholds should be chosen as
\be\label{thresh-lm1-nonunif}
\lambda(j)=\rho\frac{\ln^{\gamma}(N)}{\sqrt{\sum^{\bf r}_{l=1}\left[\frac{1}{N}+\frac{2^{-j}\pi_l}{N^{\alpha_l}}\right]^{-1}2^{-2\nu_l j}}},
\ee
where $\gamma=0.5$ under $i.i.d.$ errors and $\gamma=1$ under long-memory errors, and $\rho$ is a positive tuning parameter. In addition, under $i.i.d.$ errors or {long-memory errors} with $\Delta=0$, the finest resolution level is chosen as 
\be \label{finresol}
2^{J}=\frac{{N}^{\frac{1}{2\nu_o+1}}}{\ln(N)},
\ee
and under long memory with $\Delta\neq 0$, it is chosen as
\be \label{finresol-lm*}
2^{J}=\frac{{N}^{\frac{1}{2\nu^*+1}}}{\ln(N)}.
\ee
Notice that under $i.i.d.$ errors or under {homoskedastic} long-memory errors ($\Delta=0$) the thresholds \fr{thresh-lm1-nonunif-iid} and the resolution level \fr{finresol} {do not depend on the knowledge of long memory parameters $\alpha_l$}. In addition, \fr{finresol} depends only on the smallest smoothness parameter of the convolution kernels across all ${\bf r}$ channels. In contrast, under {heteroskedastic} long-memory errors ($\Delta\neq 0$), the choices of the thresholds \fr{thresh-lm1-nonunif} and the finest resolution level \fr{finresol-lm*} depend on the knowledge of long memory parameters $\alpha_l$. Nonetheless, estimator \fr{fest} is adaptive {with respect to the unknown smoothness parameter $s$ of the signal $f$}.
 \begin{lemma} \label{lem:Lardev-nonunif-iid}
Let $\widehat{\beta_{\omega}}$ be defined in \fr{betest1}, and let assumptions ${\bf{A.1}}$-${\bf{A.4}}$ hold. {Assume $h_l(t)$ are known}, with $l=1, 2, \cdots, {\bf r}$ and choose the thresholds as \fr{thresh-lm1-nonunif-iid} under $i.i.d.$ errors or {long-memory errors} with $\Delta=0$ and \fr{thresh-lm1-nonunif} under long-memory errors with $\Delta\neq0$. Then, as $N \rightarrow \infty$, one has 
\be  \label{Largdev-iid}
\Pr \left(| \widehat{\beta_{\omega}}-\beta_{\omega} |> 1/2\lambda({j})\right)\asymp \left[\frac{1}{N}\right]^{2\tau},
\ee
where $\tau$ is a positive constant that is independent of $N$. 
\end{lemma}
\section{Asymptotic Properties of the Wavelet Estimator}
We define the minimax $L^2$-risk over a set $\Theta$ as 
\beqns
\mathcal{R}(\Theta)=\inf_{\tilde{f}}\sup_{f\in \Theta}\EE\|\tilde{f}-f\|^2,
\eeqns
where the infimum is taken over all possible estimators $\tilde{f}$ of $f$ and
\be \label{sqnorm}
\|f_1-f_2\|^2=\int^1_0|f_1(x)-f_2(x)|^2dx.
\ee
{\bf Assumption A.5.} \label{A5} The function $f(t)$ belongs to a Besov space. In particular, if $s_{o} \geq s$, its wavelet coefficients $ \beta_{jk}$ satisfy
\be  \label{eq11}
 B^{s}_{p, q}(A)=\left \{ f \in L^2([0, 1]): \left( \sum_{j} 2^{js^{*}q}\left (\sum_{k}| \beta_{jk}|^{p}\right)^{q/{p}}\right )^{1/q} \leq A\right \}.
\ee
Define the quantities
\beqn  \label{eq10}
  p^*&=&\frac{1}{p}-\frac{1}{2},\\
 s^{*}&=&s-p^*,\\
  s^{'}&=&s+1/2 - 1/p', \ \ p'=(p\wedge 2).
 \eeqn
 According to Yang~(2001), the minimax risk of estimating a regression function under random design and long-memory errors is the sum of two terms; the minimax rate of the same class {in random design} but under $i.i.d.$ errors, and the rate of estimating the mean of the regression function under {long-memory} errors. Therefore, under simultaneous deconvolution regression with long-memory errors, we have ${\bf r}$ means to estimate and thus the lower bounds of the risk will reflect that feature. 
\begin{theorem}\label{th:lowerbds-nonunif-iid} Let $s \geq \max\{\frac{1}{p}, \frac{1}{2} \}$ with $1 \leq p,q \leq \infty$, and $A> 0$. Let  Assumptions {\bf{A.1-A.5}}  hold. {Assume $f\in B^{s}_{p, q}(A)$ and $h_l(t)$ are known}, with $l=1, 2, \cdots, {\bf r}$. Then, under $i.i.d.$ errors, as $N\rightarrow \infty$,
 \be \label{lowerbds1-nonunof-iid}
 \mathcal{R}({{B^{s}_{p, q}}}(A))\geq C\left\{ \begin{array}{ll} 
A^2 \left[\frac{\sigma^2_o}{A^2N}\right]^{\frac{2s}{2s+2\nu_o +1}}    , & \mbox{if}\ \   {s}>  p^*(2\nu_o+1), \\
A^2 \left[\frac{\sigma^2_o\ln(N)}{A^2N}\right]^{\frac{2s^*}{2s^*+2\nu_o}}   , & \mbox{if}\ \   {s}\leq p^*(2\nu_o+1).
  \end{array} \right.
   \ee
   In addition, under long memory errors, as $N\rightarrow \infty$,
   \be \label{lowerbds1-nonunif-lm}
 \mathcal{R}({{B^{s}_{p, q}}}(A))\geq C\left\{ \begin{array}{ll} 
A^2 \left[\frac{\sigma^2_{o}}{A^2N}\right]^{\frac{2s}{2s+2\nu_o+1}}\vee \frac{1}{N^{\alpha_o}}, & \mbox{if}\ \   {s}>  p^*(2\nu_o+1), \\
A^2 \left[\frac{\sigma^2_{o}\ln(N)}{A^2N}\right]^{\frac{2s^*}{2s^*+2\nu_o}}\vee \frac{1}{N^{\alpha_o}}     , & \mbox{if}\ \  {s}\leq p^*(2\nu_o+1).
 \end{array} \right.
   \ee
  \end{theorem}
  {Observe under $i.i.d.$ errors the usual behavior of one elbow effect and two different phases. The lower bounds are identical to those in De Canditiis and Pensky~(2006) under simultaneous deconvolution in regular design and $i.i.d.$ errors. In the case of ${\bf r}=1$, the bounds are identical to those in Benhaddou~(2025) under standard deconvolution in random design and $i.i.d.$ errors, and those in (Donoho~(1995a), Abramovich and Silverman~(1998), Johnstone et al.~(2004)), under standard deconvolution in regular design and $i.i.d.$ errors. In addition, under $i.i.d.$ errors the bounds depend only on the smallest degree of ill-posedness $\nu_o$ of functions $g_l$ across all ${\bf r}$ channels.} On the other hand, under long-memory errors, we see the same elbow effects observed under the $i.i.d. $ case in addition to a long-memory phase. {This latter} phase depends only on the channel that corresponds to the strongest long memory. More precisely, the lower bounds choose the worst case between the lower bounds under $i.i.d.$ and the mean squared error for parametric mean estimation under the strongest long-memory across the channels. {This behavior of the lower bounds under long-memory errors is completely new and it tells us that randomization in simultaneous deconvolution model when errors are long-memory dependent may not improve the estimation if long memory is sufficiently strong in at least one channel. This contradicts what has been encountered in other nonparametric models, such as nonparametric regression under random design and long-memory errors, (e.g., Yang~(2001), Kulik and Raimondo~(2009), and Benhaddou and Liu~(2024)), that randomization improves estimation}.
 \begin{theorem} \label{th:upperbds-2-nonunif-iid}
Suppose assumptions ${\bf{A.1}}$-${\bf{A.5}}$ hold. Assume  $f\in B^{s}_{p, q}(A)$ {and $h_l(t)$ are known}, with $l=1, 2, \cdots, {\bf r}$. Let  $\widehat{f}(t)$ be the estimator defined in \fr{fest}. Choose under $i.i.d.$ errors $\lambda(j)$ as \fr{thresh-lm1-nonunif-iid} and $J$ as \fr{finresol}. Also, choose $\rho$ in \fr{thresh-lm1-nonunif-iid} such that 
\beqn\label{rho-b}
\rho^2> 18\max\left\{\frac{M^2\max_l\left\{ K^2_{l2}\left(\frac{2\pi}{3}\right)^{-4\nu_l}\right\}}{\min_l \left\{K^2_{l1}m_{l1}\left(\frac{8\pi}{3}\right)^{-4\nu_l}\right\}}, \frac{\max_l\left\{ K^2_{l2}\left(\frac{2\pi}{3}\right)^{-4\nu_l}\right\}}{\min_l \left\{K^2_{l1}\left(\frac{8\pi}{3}\right)^{-4\nu_l}\right\}}\right\}.
\eeqn
Then, as $N \rightarrow \infty$, one has
 \be \label{upperbds-2-iid}
\EE\|\widehat{f}-f\|^2\leq CA^2\left\{ \begin{array}{ll} 
 \left[\frac{\ln(N)}{A^2N}\right]^{\frac{2s}{2s+2\nu_o +1}}  , & \mbox{if}\ \   {s}> p^*(2\nu_o+1), \\
 \left[\frac{\ln(N)}{A^2N} \right]^{\frac{2s'}{2s'+2\nu_o}}\left[\ln(N)\right]^{d}, & \mbox{if}\ \  {s}\leq p^*(2\nu_o+1),
  \end{array} \right.
\ee
where
\be\label{d-iid-unif}
d=\II\left(s=(2\nu_o+1)p^*\right). 
\ee
In addition, under long memory errors with $\Delta=0$, choose $\lambda(j)$ as \fr{thresh-lm1-nonunif-iid}  and $J$ according to \fr{finresol}. Then, if $\rho$ satisfies \fr{rho-b}, as $N \rightarrow \infty$, one has  
\be \label{upperbds-2-detl0}
 \EE\|\widehat{f}-f\|^2\leq CA^2\left\{ \begin{array}{ll} 
 \left[\frac{\ln(N)}{A^2N}\right]^{\frac{2s}{2s+2\nu_o +1}}\vee  \left[\frac{\ln(N)}{{A}^2N^{\alpha_o}}\right], & \mbox{if}\ \   p^*(2\nu_o+1)< {s}, \\
 \left[\frac{\ln(N)}{A^2N} \right]^{\frac{2s'}{2s'+2\nu_o}}\left[\ln(N)\right]^{d}\vee \left[\frac{\ln(N)}{{A}^2N^{\alpha_o}}\right], & \mbox{if}\ \  {s}\leq p^*(2\nu_o+1).
 \end{array} \right.
\ee
Furthermore, under long memory errors with $\Delta\neq 0$, choose $\lambda(j)$ as \fr{thresh-lm1-nonunif}  and $J$ according to \fr{finresol-lm*}. Then, if $\rho$ satisfies \fr{rho-b}, as $N \rightarrow \infty$, one has  
\be \label{upperbds-2}
 \EE\|\widehat{f}-f\|^2\leq CA^2\left\{ \begin{array}{ll} 
 \left[\frac{\ln^2(N)}{A^2N}\right]^{\frac{2s}{2s+2\nu^* +1}}\vee  \left[\frac{\ln^2(N)}{{A}^2N^{\alpha^*}}\right]^{\frac{2s}{2s+2\nu^*}}   , & \mbox{if}\ \   p^*(2\nu^*+1)< {s}, \\
 \left[\frac{\ln^2(N)}{A^2N} \right]^{\frac{2s'}{2s'+2\nu^*}}\left[\ln(N)\right]^{d'}\vee \left[\frac{\ln^2(N)}{{A}^2N^{\alpha^*}}\right]^{\frac{2s}{2s+2\nu^*}}, & \mbox{if}\ \  2\nu^* p^*<{s}\leq p^*(2\nu^*+1),\\
 \left[\frac{\ln^2(N)}{A^2N} \right]^{\frac{2s'}{2s'+2\nu^*}}\left[\ln(N)\right]^{d'}\vee \left[\frac{\ln^2(N)}{A^2N^{\alpha^*}} \right]^{\frac{2s'}{2s'+2\nu^*-1}}\left[\ln(N)\right]^{d''}, & \mbox{if}\ \  2\nu^* p^*\geq{s},
   \end{array} \right.
\ee
 where
\be\label{d-lm-unif}
d''=\II\left(s=2\nu^* p^*\right), \ \ \ d'=\II\left(s=(2\nu^*+1)p^*\right). 
\ee
\end{theorem}
{\bf Theorems} {\bf \ref{th:lowerbds-nonunif-iid}} and {\bf \ref{th:upperbds-2-nonunif-iid}} imply that, for the $L^2$-risk under $i.i.d.$ errors or {homoskedastic} long-memory errors ($\Delta=0$), estimator \fr{fest} with {choices} \fr{thresh-lm1-nonunif-iid} and \fr{finresol} is asymptotically near-optimal {within a logarithmic factor of $N$} over a wide range of Besov balls. Under $i.i.d.$ errors, the results agree with De Canditiis and Pensky~(2006) in their treatment of simultaneous deconvolution under regular design and $i.i.d.$ errors. For ${\bf r}=1$, these results also agree with Benhaddou~(2025) under standard deconvolution with random design, and (Donoho~(1995a), Abramovich and Silverman~(1998), Johnstone et al.~(2004)), for standard deconvolution under fixed design and $i.i.d.$ errors. In addition, the convergence rates depend only on the smallest degree of ill-posedness of the convolution kernels across the channels, a behavior that was first observed in De Canditiis and Pensky~(2006). Furthermore, under long-memory errors, the convergence rates are new and more complicated than those found in Kulik, Sapatinas and Wishart~(2015) {in their treatment of simultaneous deconvolution under regular design and long-memory errors}. In particular, when $\Delta=0$, the convergence rates depend only on the smallest long-memory parameter, $\alpha_o$ (the strongest long memory) and the lowest degree of ill-posedness, $\nu_o$, across the channels. {This tells us that when long memory is relatively weak across the channels, the convergence rates coincide with De Canditiis and Pensky~(2006), in their treatment of the same model under regular design and $i.i.d.$ errors}. {However, when long memory is sufficiently strong in at least one channel, the convergence rates depend only on the strongest long memory}. {This tells us that when long memory is sufficiently strong in at least one channel, the convergence rates are completely different from De Canditiis and Pensky~(2006), in which we see that the convergence rates depend on regularity of $f$ and the smallest degree of ill-posedness across the channels, and different from Kulik, Sapatinas and Wishart~(2015) in which the convergence rates depend on a delicate balance between the regularity of $f$, the degrees of ill-posedness of the convolutions and the long-memory parameters across the channels}. In addition, when $\Delta\neq 0$,  the {upper bounds of the $L^2$-risk} depend on the long memory parameter and the smoothness $\nu^*$ associated with the channel that corresponds to choice \fr{l-optimal*}. At each elbow, when the long memory is {sufficiently} strong, which corresponds to lower values of $\alpha^*$ relative to the other parameters, the {upper bounds} depend on $\alpha^*$ in addition to $s$ and $\nu^*$. In particular, when $p^*(2\nu^*+1)< {s}$, the {bounds} depend on $\alpha^*$ whenever $\alpha^*< \frac{2s+2\nu^*}{2s+2\nu^*+1}$, whereas, when $s\leq 2\nu^* p^*$, { they} depend on $\alpha^*$ whenever $\alpha^*< \frac{2s'_1+2\nu^*-1}{2s'+2\nu^*}$. {These upper bounds} are different from those found in Kulik, Sapatinas and Wishart~(2015) {in which the bounds depend on long memory regardless of its strength}. Notice that {in case $\Delta\neq 0$, the upper bounds} in \fr{upperbds-2} are not optimal {since they do not match with the lower bounds}. \\
In Kulik, Sapatinas and Wishart~(2015), the upper bounds depend on the long memory parameters  and the convolution kernels' smoothness parameters according to the rule
\beqn\label{l-kulsapwis}
  l_{*}=\arg \min_{1\leq l\leq {\bf r}}\left\{\frac{2^{2\nu_l+\alpha_l}}{N^{\alpha_l}}\right\}.
\eeqn
This rule is {slightly} different from {rule} \fr{l-optimal*} {under long-memory errors with} $\Delta \neq 0$, but it is very close {and should agree on the choice of the best channel}. 
 \begin{remark}
  {In practice, the design densities $h_l$, $l=1, 2,\cdots, {\bf r}$, may not be known in advance, but can be estimated from the data using some consistent estimators $\widehat{h}_l$, and replacing $h_l$ by estimates $\widehat{h}_l$ in \fr{betest1}. Provided that $h_l$ are sufficiently smooth, their estimation will not substantially affect thresholds \fr{thresh-lm1-nonunif-iid} and \fr{thresh-lm1-nonunif}, and therefore, it will not negatively affect the convergence rates of estimating $f$. However, if $h_l$ are not sufficiently smooth, thresholds \fr{thresh-lm1-nonunif-iid} and \fr{thresh-lm1-nonunif} will be considerably affected as their choice will depend on the regularities of ${h}_l$, and this subsequently will have a detrimental effect on the convergence rates of estimating $f$.}
   \end{remark}	  
   \begin{remark}	  
{ In practice, the parameters $\alpha_l$ may not be known in advance, but can be {consistently} estimated. {In case $\alpha_l$ are unknown, then their estimates $\widehat{\alpha}_l$ will be used in the weights \fr{theta1-lm-unif} and thresholds \fr{thresh-lm1-nonunif}  instead of $\alpha_l$. It is worth noting that this will not affect the convergence rates of estimating $f$ since consistent estimators $\widehat{\alpha}_l$ converge to $\alpha_l$ at parametric rates, the quantity }
\beqns
\frac{1}{N^{\widehat{\alpha}_l}}=\frac{1}{N^{\widehat{\alpha}_l}}\frac{N^{\alpha_l}}{N^{\alpha_l}}=\frac{1}{N^{{\alpha}_l}}\frac{N^{\alpha_l}}{N^{\widehat{\alpha}_l}}=\frac{1}{N^{{\alpha}_l}}O(1), \ \ \ as\ \ \  N\rightarrow \infty. 
\eeqns}
 \end{remark}
 \section{Discussion }
   \noindent  In this paper, we constructed wavelet hard-thresholding estimators for simultaneous deconvolution model under random design in which the design densities differ across the channels. We derived lower and upper bounds for the $L^2$-risk under the two noise scenarios, and showed that the estimator is asymptotically near-optimal {up to a logarithmic factor of $N$} in the minimax sense under $i.i.d.$ errors and {homoskedastic} long-memory errors ($\Delta=0$). However, under {heteroskedastic} long-memory errors ($\Delta\neq0$), the estimator is {sub-optimal}. {The study of simultaneous deconvolution under random design is motivated by multi-LiDAR applications in which raw data points are inherently irregular.} \\
  \noindent  {When $\nu_l=\nu$, $l=1, 2, \cdots, {\bf r}$, under long-memory errors with $\Delta=0$, the minimax rates are dominated by the channel that corresponds to the strongest long memory, the smallest $\alpha_l$. In contrast, under long-memory errors with $\Delta\neq 0$,  the convergence rates are dominated by the channel that corresponds to the weakest long memory, the largest $\alpha_l$. On the other hand, when $\alpha_l=\alpha$, $l=1, 2, \cdots, {\bf r}$, the rates are dominated by the channel that corresponds to the smallest $\nu_l$, $l=1, 2, \cdots, {\bf r}$. For the general case, under long-memory errors with $\Delta=0$, the convergence rates are dominated by the channel that corresponds to the lowest $\nu_l$ and the channel that corresponds to the smallest $\alpha_l$. However, under long-memory errors with $\Delta\neq 0$, the rates may be dominated by the channel corresponding to neither the weakest long memory $\alpha_l$ nor the smallest $\nu_l$, $l=1, 2, \cdots, {\bf r}$.}  \\
1. De Canditiis and Pensky~(2006) consider simultaneous deconvolution under regular, evenly-spaced design and  $i.i.d.$ errors to construct a wavelet hard-thresholding estimator and show that it is asymptotically near optimal in the minimax sense. Their convergence rates depend only on the smallest smoothness $\nu_l$ across the channels, so our results under $i.i.d.$ errors are identical to theirs. Under long-memory errors with $\Delta=0$, when {long-memory is relatively weak across the channels}, our convergence rates are identical to theirs. {However, when long memory is sufficiently strong in at least one channel, the convergence rates are completely different from theirs, and depend only on the strongest long memory across the channels}. On the other hand, under long-memory errors with $\Delta\neq 0$, even when the parameter $\alpha^*$ in choice \fr{l-optimal*} is {relatively} large,  our convergence rates may still not be identical to theirs, considering that the smoothness $\nu^*$ in \fr{l-optimal*} may not be the smallest {across the channels}. \\
2.  For ${\bf r}=1$, our convergence rates {are identical}  to those in Benhaddou~(2025), for the treatment of standard deconvolution under random design, and {both $i.i.d.$ errors and long-memory errors with $\Delta\neq 0$}. {It is worth noting that Benhaddou~(2025) did not consider case $\Delta=0$ under long-memory errors, and their lower bounds are not as tight as those in {\bf Theorems} {\bf \ref{th:lowerbds-nonunif-iid}} and could be improved to match ours}. In addition, with $\nu=0$ and {long-memory errors}, our rates are identical to those observed in nonparametric wavelet regression, {such as} Yang~(2001), Kulik and Raimondo~(2009), and Benhaddou and Liu~(2024).\\
3.  Kulik, Sapatinas and Wishart~(2015) {investigate} simultaneous deconvolution under regular, evenly-spaced design and long memory errors. Their rates depend on the channel that corresponds to rule \fr{l-kulsapwis}, {which somewhat agrees with our choice \fr{l-optimal*}} under long-memory errors with $\Delta\neq 0$. {In homoskedastic long-memory ($\Delta=0$) case, when long-memory is sufficiently strong in at least one channel, our convergence rates depend only on the long-memory parameters. In addition, when long memory is relatively weak across the channels, our convergence rates don't depend on long memory at all, but do depend on the regularity of $f$ and the degrees of ill-posedness of the convolution operators, while in their case the rates depend on long-memory parameters in addition to the regularity of $f$ and the degrees of ill-posedness of the convolution operators}. {In heteroskedastic long-memory ($\Delta\neq 0$), when long-memory parameter for channel corresponding to choice \fr{l-optimal*} is relatively large (weaker long memory), our convergence rates don't depend on long memory at all, but do depend on the other parameters, regularities of $f$ and $g_{l^*}$. However, when long-memory parameter $\alpha^*$ is relatively small (stronger long memory), the rates depend on $\alpha^*$ in addition to $s$ and $\nu^*$, but in a way that is different from  Kulik, Sapatinas and Wishart~(2015). Our convergence rates are faster than theirs in this case. }
 \section{Proofs }
In order to prove {\bf Theorems \ref{th:lowerbds-nonunif-iid}}, we use the following lemma.  
\begin{lemma} (Lemma $A.1$ of Bunea et al.~\cite{bun}) 
\label{lem:Bunea} 
Let $\Te$ be a set of functions of cardinality $\card(\Te)\geq 2$ such that\\
(i) $\|f-g\|_2^2 \geq 4\delta^2, \ for\  f, g \in \Te, \ f \neq g, $\\
(ii) the Kullback divergences $\mathcal{K}(P_f, P_g)$ between the measures $P_f$ and $P_g$ 
satisfy the inequality $\mathcal{K}(P_f, P_g) \leq \log(\card(\Te))/16,\ for\ f,\ g \in \Te$.\\
Then, for some absolute positive constant $C_5$, one has 
$$
 \inf_{f_n}\sup_{f\in \Te} \EE_f\|f_n-f\|_2^2 \geq C_5 \delta^2,
$$
where $\inf_{f_n}$ denotes the infimum over all estimators.
\end{lemma}
 {\bf Proof of Theorem \ref{th:lowerbds-nonunif-iid}}. {\bf The $i.i.d.$ errors case}. \\
 {\bf Dense Case}: Using the same test functions $f_{\omega}$ and $f_{\tilde{\omega}}$ as in the proof of {\bf Theorem 3.1} of Benhaddou~(2016) under the dense case, it can be shown that the $L^2$-norm of the difference satisfies 
\be\label{l2diff}
\| f_{\omega}-f_{\tilde{\omega}}\|^2 \geq  \frac{2^{j}\rho^2_{j}}{8}.
\ee
In order to apply Lemma \ref{lem:Bunea}, one needs to verify condition (ii). Denote, for $l=1, 2, \cdots, {\bf r}$,  
\be
q_{l\omega}(t_{li})=\int^1_0f_{\omega}(s)g_l(t_{li}-s)ds, \ i=1, 2, \cdots, N.
\ee
\be
q_{l\tilde{\omega}}(t_{li})=\int^1_0f_{\tilde{\omega}}(s)g_l(t_{li}-s)ds, \ i=1, 2, \cdots, N.
\ee
Define the quantities $\pmb{h}_{li}= q_{l\omega}(t_{li})+ \sigma_l(t_{li})\varepsilon^{(l)}_{i}$, $i=1, 2, \cdots, N$, where $\{\varepsilon^{(l)}_{i}\}_{(i)\geq 1}$ is $i.i.d.$ Gaussian sequence. Let $P_{f_{\omega}}$ be the probability law of the process $\pmb{h}_{i}=\left(\pmb{h}_{1i}, \pmb{h}_{2i}, \cdots, \pmb{h}_{{\bf r}i}\right)$ under the hypothesis $f_{\omega}$. Then, by Assumptions ${\bf A.2}$-${\bf A.5}$, the Kullback divergence can be written as 
\beqn
\mathcal{K}(P_{f_{\omega}}, P_{f_{\tilde{\omega}}})&\leq& \frac{1}{2}\left\{N\sum^{\bf r}_{l=1}\EE\left[\frac{q_{l\omega}(t_1)-q_{l\tilde{\omega}}(t_1)}{\sigma_l^2(t_1)}\right]_{h_l}^2\right\}\nonumber\\
&\leq&\sum^{\bf r}_{l=1} \frac{m_{l2}}{2}N\|\sigma^{-2}_{l}(t)\|_{\infty}\rho^2_{j}2^{j}\left(\frac{2\pi}{3}\right)^{-2\nu_l}2^{-2\nu_l j}\nonumber\\
&\leq & \sum^{\bf r}_{l=1} \frac{m_{l2}}{2}N\sigma^{-2}_{lo}\rho^2_{j}2^{j}\left(\frac{2\pi}{3}\right)^{-2\nu_l}2^{-2\nu_l j}\\
&=&CN2^j\rho^2_{j}\sum^{\bf r}_{l=1}\sigma^{-2}_{lo}2^{-2\nu_l j}.
\eeqn
Now, to apply {\bf Lemma \ref{lem:Bunea}}, optimally choose 
\be \label{klblam}
A^2 2^{-j(2s+2\nu_o+1)}2^{j}\frac{N}{\sigma_o^2}\leq \pi_o2^{j},
\ee
where $\sigma_o$ is the ${\sigma_{lo}}$ that corresponds to the equation $\pmb{h}_{li}$ for which $\nu_o=\min\{\nu_1, \nu_2,\cdots, \nu_{\bf r}\}$, that is,
\be \label{Lop1}
2^{j} = C\left[\frac{A^2N}{\sigma_o^2}\right]^{\frac{1}{2s+2\nu_o +1}}.
\ee
Hence, the lower bound is
\be \label{cas1-lm1-iid}
\delta^2=C 
A^2 \left[\frac{\sigma^2_o}{A^2N} \right]^{\frac{2s}{2s+2\nu_o +1}}.
\ee
 2: {\bf Sparse Case}. We follow the same procedure as in the {\bf Dense Case} and use the test functions from Benhaddou~(2016) to obtain the lower bounds
\be \label{cas2-lm1-iid}
\delta^2=CA^2 \left[\frac{\sigma^2_o}{A^2N} \right]^{\frac{2s^*}{2s^*+2\nu_o}}.
\ee
Now, we choose the highest between \fr{cas1-lm1-iid} and \fr{cas2-lm1-iid}.  
 $\Box$\\
 {\bf Long memory errors case}. 
Let, for $l=1, 2, \cdots, {\bf r}$, $\Sigma^{(l)}_N$ denote the $N\times N$ covariance matrix of errors $\varepsilon^{(l)}_{i}$, and $\left[\xi^{(l)}\right]^{-1}_{ii'}$ be the $ii'$ element of the inverse of $\Sigma^{(l)}_N$. Also, let ${\bf1}$ denote the N-dimensional vector whose elements are all equal to 1. \\
  {\bf Dense Case}: The $L^2$-norm of the difference satisfies \fr{l2diff}. Consider the quantities $\pmb{h}_{li}$ defined above, but now assume $\{\varepsilon^{(l)}_{i}\}_{(i)\geq 1}$ is stationary Gaussian sequence that satisfies property \fr{longm1}, and such that $\EE\left[\varepsilon^{(l)}_{i}\varepsilon^{(l')}_{i}\right]=0$ for $l\neq l'$. Let $P_{f_{\omega}}$ be the probability law of the process $\pmb{h}_{i}=\left(\pmb{h}_{1i}, \pmb{h}_{2i}, \cdots, \pmb{h}_{{\bf r}i}\right)$ under the hypothesis $f_{\omega}$. Then, by Assumptions ${\bf A.2}$-${\bf A.5}$, and Holder's inequality, the Kullback divergence can be written as 
\beqn
\mathcal{K}(P_{f_{\omega}}, P_{f_{\tilde{\omega}}})&\leq& \frac{1}{2}\sum^{\bf r}_{l=1}\left\{\EE\left[\frac{q_{l\omega}(t_{li})-q_{\tilde{\omega}}(t_{li})}{\sigma_l^2(t_{li})}\right]^2_{h_l}\sum^N_{i=1}\left[\xi^{(l)}\right]^{-1}_{ii}+\left[\EE\left(\frac{q_{l\omega}(t_{li})-q_{l\tilde{\omega}}(t_{li})}{\sigma_l(t_{li})}\right)_{h_l}\right]^2\sum^N_{i\neq i'}\left[\xi^{(l)}\right]^{-1}_{ii'} \right\}\nonumber\\
&\leq& \sum^{\bf r}_{l=1}\frac{m_{2}}{2}\|\sigma^{-2}_{l}\|_{\infty}\rho^2_{j}2^{j}\left(\frac{2\pi}{3}\right)^{-2\nu_l}2^{-2\nu_l j}\left\{ Trace\left(\left[\Sigma^{(l)}_N\right]^{-1}\right)+{\bf1}^T\left[\Sigma^{(l)}_N\right]^{-1}{\bf1}\right\}\nonumber\\
&\leq &  \sum^{\bf r}_{l=1}\frac{m_{2}}{2}\sigma^{-2}_{lo}\rho^2_{j}2^{j}\left(\frac{2\pi}{3}\right)^{-2\nu_l}2^{-2\nu_l j}\left\{ N+\pi^{-1}_{l}N^{\alpha_l} \right\}.
\eeqn
To see that ${\bf1}^T\left[\Sigma^{(l)}_N\right]^{-1}{\bf1}\sim\pi^{-1}_{l}N^{\alpha_l}$, note that processes satisfying condition \fr{longm1} also have the property
\be\label{lambda-sigrma}
C_{l1}N^{1-\alpha_l}\leq \lambda_{\min}\left[\Sigma^{(l)}_N\right]\leq \lambda_{\max}\left[\Sigma^{(l)}_N\right]\leq C_{l2}N^{1-\alpha_l},
\ee
where $C_{l1}$ and $C_{l2}$ are positive constants and $\lambda_{\min}\left[\Sigma^{(l)}_N\right]$ and $\lambda_{\max}\left[\Sigma^{(l)}_N\right]$ are the smallest and the largest eigenvalue of $\Sigma^{(l)}_N$, respectively (e.g., see Benhaddou et al.~(2014), section 2). Therefore, by \fr{lambda-sigrma}
\beqns
{N}{C_{l1}}N^{1-\alpha_l}\leq {\bf1}^T\left[\Sigma^{(l)}_N\right]{\bf1}&\leq&{N}{C_{l2}}N^{1-\alpha_l},\nonumber\\
\frac{N}{C_{l2}}N^{\alpha_l-1}\leq {\bf1}^T\left[\Sigma^{(l)}_N\right]^{-1}{\bf1}&\leq&\frac{N}{C_{l1}}N^{\alpha_l-1}.
\eeqns
Or
\beqns
{\bf1}^T\left[\Sigma^{(l)}_N\right]{\bf1}&\sim&\pi_{l}N^{2-\alpha_l}, \nonumber\\
{\bf1}^T\left[\Sigma^{(l)}_N\right]^{-1}{\bf1}&\sim&\pi^{-1}_{l}N^{\alpha_l}. 
\eeqns
Now, to apply {\bf Lemma \ref{lem:Bunea}}, optimally choose 
\be \label{klblam}
A^2\sum^{\bf r}_{l=1} 2^{-j(2s+2\nu_l+1)}2^{j}\frac{N}{\sigma_{lo}^2}+ A^2 \sum^{\bf r}_{l=1} 2^{-j(2s+2\nu_l+1)}2^{j}\frac{N^{\alpha_l}}{\pi_{l}\sigma_{lo}^2}\leq \pi_o2^{j},
\ee
that is, by choosing the channel that satisfies \fr{l-optimal*}, the optimal level $j$ is 
\be \label{Lop1}
2^{j} = C\left[\frac{A^2N}{\sigma_{o}^2}\right]^{\frac{1}{2s+2\nu_o +1}}.
\ee
Hence, the lower bounds are
\be \label{cas1-lm1}
\delta^2=C \left\{ 
A^2 \left[\frac{\sigma^2_{o}}{A^2N} \right]^{\frac{2s}{2s+2\nu_o +1}} \vee \max_{l\leq {\bf r}}\left\{\frac{1}{N^{\alpha_l}}\right\} \right\}.
\ee
 2: {\bf Sparse Case}. We follow the same procedure as in the {\bf Dense Case} and use the test functions from Benhaddou~(2016) to obtain the lower bounds
\be \label{cas2-lm1}
\delta^2=C \left\{ 
A^2 \left[\frac{\sigma^2_{o}}{A^2N} \right]^{\frac{2s^*}{2s^*+2\nu_o}} \vee \max_{l\leq {\bf r}}\left\{\frac{1}{N^{\alpha_l}}\right\}\right\}.
\ee
Now, we choose the highest between \fr{cas1-lm1} and \fr{cas2-lm1}.  
 $\Box$\\
 {\bf Proof of Lemma \ref{lem:Var-nonunif}.} Notice that with $\widehat{\beta^*_{\omega}}$  defined in \fr{betest1}, one has 
\beqn \label{alef-iid}
\widehat{\beta^*_{\omega}}-\beta^*_{\omega}&=& \frac{1}{N}\sum^{\bf r}_{l=1}\sum^N_{i=1}\left[\frac{q_l(t_{li})U^{(l)}_{\omega}(t_{li})}{h_l(t_{li})}-\EE\left[\frac{q_l(t_l)U^{(l)}_{\omega}(t_l)}{h_l(t_{l})}\right]\right]\nonumber\\
 &+& \frac{1}{N}\sum^{\bf r}_{l=1}\sum^N_{i=1}\frac{U^{(l)}_{\omega}(t_{li})\sigma_l(t_{li})}{h_l(t_{li})}\varepsilon^{(l)}_{i}=\Delta_{1}+\Delta_2.
 \eeqn
Under $i.i.d.$ errors, both $\Delta_{1}$ and $\Delta_{2}$ are the sums of zero-mean independent and identically distributed components. Therefore, taking expectation of the square of $\Delta_{1}$, using Cauchy-Schwarz inequality, the set \fr{w-set} and condition \fr{gm},  yields
\beqn
\EE\left[\Delta^2_{1}\right]
&\leq& \frac{1}{N}\|f\|^2_{\infty}\max_{ l\leq {\bf r}}\|g_l\|^2_{\infty}\sum^{\bf r}_{l=1}\int^1_0\left|\frac{U^{(l)}_{\omega}(t)}{h_l(t)}\right|^2h_l(t)dt\nonumber\\
&\leq & \frac{1}{N}\|f\|^2_{\infty}\max_{ l\leq {\bf r}}\|g_l\|^2_{\infty}\sum_{m\in W_{j}}\sum^{\bf r}_{l=1}\frac{1}{m_{l1}}\frac{\theta^2_l\left|\tilde{g}_l(m)\right|^2\left|\tilde{\psi}_{\omega}(m)\right|^2}{\left[\sum^{\bf r}_{l=1}\theta_l\left|\tilde{g}_l(m)\right|^2\right]^2}\nonumber\\
&\leq&\frac{C}{N}\frac{\sum^{\bf r}_{l=1}m^{-1}_{l1}\theta^2_l2^{-2\nu_lj}}{\left[\sum^{\bf r}_{l=1}\theta_l2^{-2\nu_lj}\right]^2}. \label{delta1-iid}
\eeqn
Also, it can be easily shown that
\beqn \label{D14}
\EE\left[\Delta^4_{1}\right]&\leq&\frac{C}{N^2}\left(\frac{\sum^{\bf r}_{l=1}m^{-1}_{l1}\theta^2_l2^{-2\nu_lj}}{\left[\sum^{\bf r}_{l=1}\theta_l2^{-2\nu_lj}\right]^2}\right)^2.
\eeqn
As for $\Delta_{2}$, since $\varepsilon^{(l)}_i$ are $i.i.d.$ mean-zero Gaussians and that $\EE\left[\varepsilon^{(l)}_{i}\varepsilon^{(l')}_{i}\right]=0$ for $l\neq l'$, the square expectation yields 
\beqn
\EE\left[\Delta^2_{2}\right]&=& \frac{1}{N}\sum^{\bf r}_{l=1}\sigma^2_{l1}\int^1_0\left|U^{(l)}_{\omega}(t)\right|^2dt\EE\left[\varepsilon^{(l)}_{1}\right]^2 \nonumber\\
&\leq & \frac{1}{N}\sum_{m\in W_{j}}\sum^{\bf r}_{l=1}\sigma^2_{l1}\frac{\theta^2_l\left|\tilde{g}_l(m)\right|^2\left|\tilde{\psi}_{\omega}(m)\right|^2}{\left[\sum^{\bf r}_{l=1}\theta_l\left|\tilde{g}_l(m)\right|^2\right]^2}\EE\left[\varepsilon^{(l)}_{1}\right]^2\nonumber\\
&\leq& \frac{C}{N}\frac{\sum^{\bf r}_{l=1}\theta^2_l2^{-2\nu_lj}}{\left[\sum^{\bf r}_{l=1}\theta_l2^{-2\nu_lj}\right]^2}, \label{delta2-iid}
\eeqn
and 
\beqn \label{D24}
\EE\left[\Delta^4_{2}\right] &=&\frac{1}{N^3}\sum^{\bf r}_{l=1}\sigma^4_{l1}\int^1_0\left|U^{(l)}_{\omega}(t)\right|^4dt\EE\left[\varepsilon^{(l)}_{1}\right]^4+\frac{1}{N^2}\sum^{\bf r}_{l=1}\left(\int^1_0\left|U^{(l)}_{\omega}(t)\right|^2dt\EE\left[\varepsilon^{(l)}_{1}\right]^2 \right)^2\nonumber\\
&\leq& \frac{C}{N^2}\sum^{\bf r}_{l=1}\left(\sum_{m\in W_{j}} \frac{\theta^2_l\left|\tilde{g}_l(m)\right|^2\left|\tilde{\psi}_{\omega}(m)\right|^2}{\left[\sum^{\bf r}_{l=1}\theta_l\left|\tilde{g}_l(m)\right|^2\right]^2}\right)^2\nonumber\\
&\leq& \frac{C}{N^2}\frac{\sum^{\bf r}_{l=1}\theta^4_l2^{-4\nu_lj}}{\left[\sum^{\bf r}_{l=1}\theta_l2^{-2\nu_lj}\right]^4}.
\eeqn
Hence,  taking $\theta_l$ according to \fr{theta1-iid}, combining \fr{delta1-iid} and \fr{delta2-iid} completes the proof of \fr{var-bet-nonunif-iid}, and combining \fr{D14} and \fr{D24} completes the proof of \fr{var2-bet-nonunif-iid}. $\Box$\\
Now, under long-memory errors, $\Delta_{1}$ does not depend on $\varepsilon^{(l)}$'s, so it will be the same as before and therefore results \fr{delta1-iid} and \fr{delta2-iid} are still valid. $\Delta_{2}$ on the other hand is different, and can be decomposed as follows 
\beqn \label{alef}
\Delta_2&=& \Delta_{21}+\Delta_{22} \nonumber\\
 &=& martingale/ long\ memory\ term.
 \eeqn
 To evaluate  $\Delta_{2}$ in this case, define for a fixed $l$, $\mathcal{F}^{(l)}_i=\sigma\left( \eta^{(l)}_i, t^{(l)}_i, \eta^{(l)}_{i-1}, t^{(l)}_{i-1}, \cdots \right)$ and let $\varepsilon^{(l)}_{i, i-1}=\varepsilon^{(l)}_i-\eta^{(l)}_i$. Note that $\varepsilon^{(l)}_{i, i-1}$ is $\mathcal{F}^{(l)}_{i-1}$-measurable and $(\eta^{(l)}_i, t^{(l)}_i)$ is independent of $\mathcal{F}^{(l)}_{i-1}$. Therefore, 
\beqn
\EE\left[ \frac{{U^{(l)}_{\omega}(t^{(l)}_i)\sigma_l(t^{(l)}_i)}}{h_l(t^{(l)}_i)}\varepsilon^{(l)}_{i}\mid \mathcal{F}^{(l)}_{i-1}\right]&=&\EE\left[\frac{U^{(l)}_{\omega}(t_l)\sigma_l(t_l)}{h_l(t_l)}\right] \EE\left[ \varepsilon^{(l)}_{i}\mid \mathcal{F}^{(l)}_{i-1}\right]\nonumber\\
&=&\EE\left[\frac{U^{(l)}_{\omega}(t_l)\sigma_l(t_l)}{h_l(t_l)}\right]  \varepsilon^{(l)}_{i, i-1}.
\eeqn
Thus, $\Delta_2$ can be partitioned as follows
\beqn
\Delta_2&=&\frac{1}{N}\sum^N_{i=1}\sum^{\bf r}_{l=1}\left[\frac{{U^{(l)}_{\omega}(t^{(l)}_i)\sigma_l(t^{(l)}_i)}}{h_l(t^{(l)}_i)}\varepsilon^{(l)}_{i}-\EE\left[ \frac{{U^{(l)}_{\omega}(t^{(l)}_i)\sigma_l(t^{(l)}_i)}}{h_l(t^{(l)}_i)}\varepsilon^{(l)}_{i}\mid \mathcal{F}^{(l)}_{i-1}\right]\right]\nonumber\\
&+& \frac{1}{N}\sum^{\bf r}_{l=1}\EE\left[\frac{U^{(l)}_{\omega}(t_l)\sigma_l(t_l)}{h_l(t_l)}\right] \sum^N_{i=1}\varepsilon^{(l)}_{i, i-1}\nonumber\\
&=&\Delta_{21}+\Delta_{22}. 
\label{delta2}
\eeqn
To evaluate $\Delta_{21}$, notice that, 
\beqns
\Delta^{(l)}_{21}=\frac{1}{N}\sum^N_{i=1}d^{(l)}_i,
\eeqns
and $\sum^N_{i=1}d^{(l)}_i$ is a martingale with 
\beqn
d^{(l)}_i&=&\frac{{U^{(l)}_{\omega}(t^{(l)}_i)\sigma_l(t^{(l)}_i)}}{h_l(t^{(l)}_i)}\varepsilon^{(l)}_{i}-\EE\left[ \frac{{U^{(l)}_{\omega}(t^{(l)}_i)\sigma_l(t^{(l)}_i)}}{h_l(t^{(l)}_i)}\varepsilon^{(l)}_{i}\mid \mathcal{F}^{(l)}_{i-1}\right]\nonumber\\
&=&\left[\frac{{U^{(l)}_{\omega}(t^{(l)}_i)\sigma_l(t^{(l)}_i)}}{h_l(t^{(l)}_i)}-\EE\left[\frac{U^{(l)}_{\omega}(t_l)\sigma_l(t_l)}{h_l(t_l)}\right]  \right]\varepsilon^{(l)}_{i, i-1}+\frac{{U^{(l)}_{\omega}(t^{(l)}_i)\sigma_l(t^{(l)}_i)}}{h_l(t^{(l)}_i)}\eta^{(l)}_i. 
\label{di1}
\eeqn
Therefore, taking the $p^{th}$ moment of \fr{di1}, for $p\geq 2$, yields
\beqn
\EE\left|d^{(l)}_i\right|^p\leq 2^{p-1}\left(\EE\left|\frac{U^{(l)}_{\omega}(t_l)\sigma_l(t_l)}{h_l(t_l)}\right|^p\left\{\EE\left|\varepsilon^{(l)}_{i, i-1}\right|^p+\EE\left|\eta^{(l)}_i\right|^p\right\}\right). \label{edp}
\eeqn
In addition,
\beqn
\EE\left[\left(d^{(l)}_i\right)^2\mid \mathcal{F}^{(l)}_{i-1}\right]=\left(\EE\left[\frac{U^{(l)}_{\omega}(t_l)\sigma_l(t_l)}{h_l(t_l)}\right]^2\left\{\left[\varepsilon^{(l)}_{i, i-1}\right]^2+\EE\left[\eta^{(l)}_1\right]^2\right\}\right). \label{edcond}
\eeqn
Therefore, by Rosenthal's inequality for martingales (with $p=2$), results \fr{edp} and \fr{edcond} and {\bf Assumption A.3}, one has
\beqns
\EE\left[ \Delta^{(l)}_{21}\right]^2&\leq& \frac{1}{N^2}\left(\EE\left[\sum^N_{i=1}\EE\left[\left(d^{(l)}_i\right)^2\mid \mathcal{F}^{(l)}_{i-1}\right]\right]+\sum^N_{i=1}\EE|d^{(l)}_i|^2\right)\nonumber\\
&\leq& \frac{3\sigma_{l1}^2}{N}\sum_{m\in W_{j}}\frac{m^{-1}_{l1}\theta^2_l\left|\tilde{g}_l(m)\right|^2\left|\tilde{\psi}_{\omega}(m)\right|^2}{\left[\sum^{\bf r}_{l=1}\theta_l\left|\tilde{g}_l(m)\right|^2\right]^2}\left[\EE|\varepsilon^{(l)}_{1, 0}|^2+\EE\left[\eta^{(l)}_1\right]^2\right].
\eeqns
Therefore, since $\varepsilon^{(1)}$, $\varepsilon^{(2)}$, $\cdots$, $\varepsilon^{({\bf r})}$ are independent, we have 
\be \label{delt2}
\EE\left[ \Delta_{21}\right]^2\leq \sum^{\bf r}_{l=1}\frac{3\sigma_{l1}^2}{N}\sum_{m\in W_{j}}\frac{m^{-1}_{l1}\theta^2_l\left|\tilde{g}_l(m)\right|^2\left|\tilde{\psi}_{\omega}(m)\right|^2}{\left[\sum^{\bf r}_{l=1}\theta_l\left|\tilde{g}_l(m)\right|^2\right]^2}\left[\EE|\varepsilon^{(l)}_{1, 0}|^2+\EE\left[\eta^{(l)}_1\right]^2\right]. 
\ee
Now, for $p=4$, one has
\beqns
\EE\left[ \Delta^{(l)}_{21}\right]^4&\leq& \frac{1}{N^4}\left(\EE\left[\sum^N_{i=1}\EE\left[\left(d^{(l)}_i\right)^2\mid \mathcal{F}^{(l)}_{i-1}\right]\right]^2+\sum^N_{i=1}\EE|d^{(l)}_i|^4\right)\nonumber\\
&\leq& \frac{\sigma_{l1}^4}{N^4}\left\{\sum_{m\in W_{j}}\frac{m^{-1}_{l1}\theta^2_l\left|\tilde{g}_l(m)\right|^2\left|\tilde{\psi}_{\omega}(m)\right|^2}{\left[\sum^{\bf r}_{l=1}\theta_l\left|\tilde{g}_l(m)\right|^2\right]^2}\right\}^2\EE\left[\sum^N_{i=1}\left[|\varepsilon^{(l)}_{1, 0}|^2+\EE\left[\eta_1\right]^2\right]\right]^2\nonumber\\
&+&\frac{8\sigma_{l1}^4}{N^4}\left\{\sum_{m\in W_{j}}\frac{m^{-1}_{l1}\theta^2_l\left|\tilde{g}_l(m)\right|^2\left|\tilde{\psi}_{\omega}(m)\right|^2}{\left[\sum^{\bf r}_{l=1}\theta_l\left|\tilde{g}_l(m)\right|^2\right]^2}\right\}^2\sum^N_{i=1}\left[\EE|\varepsilon_{1, 0}|^4+\EE\left[\eta_1\right]^4\right].
\eeqns
Therefore, it can be shown that
\be \label{delt24}
\EE\left[ \Delta_{21}\right]^4\leq \frac{C}{N^2} \sum^{\bf r}_{l=1}\left\{\sum_{m\in W_{j}}\frac{m^{-1}_{l1}\theta^2_l\left|\tilde{g}_l(m)\right|^2\left|\tilde{\psi}_{\omega}(m)\right|^2}{\left[\sum^{\bf r}_{l=1}\theta_l\left|\tilde{g}_l(m)\right|^2\right]^2}\right\}^2. 
\ee
Now, for the long-memory term, notice that under property $\Delta=0$, $\Delta_{22}$ in \fr{delta2} is equal to zero
\be\label{del22-hom}
\Delta_{22}=0.
\ee
Hence, choosing $\theta_l$ according to  \fr{theta1-iid}  and combining \fr{delta1-iid}, \fr{delta2-iid} and \fr{del22-hom}, completes the proof of \fr{var-bet-nonunif-iid}, and combining \fr{D14}, \fr{D24} and \fr{del22-hom} completes the proof of \fr{var2-bet-nonunif-iid}. However, under condition $\Delta\neq0$ and therefore by{\bf Assumption A.3}, one has
\beqn \label{delt'22}
\EE\left[\Delta_{22}\right]^2&=&\frac{1}{N^2}\sum^{\bf r}_{l=1}\left(\EE\left[\frac{U^{(l)}_{\omega}(t_l)\sigma_l(t_l)}{h_l(t_l)}\right]\right)^2\Var\left[\sum^N_{i=1}\varepsilon^{(l)}_{i, i-1}\right]\nonumber\\
&\leq &\sum^{\bf r}_{l=1} \left(\sum_{m\in W_{j}}\frac{\theta_l\left|\tilde{g}_l(m)\right|\tilde{\psi}_{\omega}(m)\overline{\tilde{\sigma}_l(m)}}{\left[\sum^{\bf r}_{l=1}\theta_l\left|\tilde{g}_l(m)\right|^2\right]}\right)^2\pi_{l} N^{-\alpha_l}\nonumber\\
&\leq&\sum^{\bf r}_{l=1}\frac{\theta^2_l\pi_{l}}{N^{\alpha_l}} \max_{m\in W_j}\left\{\frac{\left|\tilde{g}_l(m)\right|^2}{\left[\sum^{\bf r}_{l=1}\theta_l\left|\tilde{g}_l(m)\right|^2\right]^2} \right\}\left(\sum_{m\in W_{j}}\tilde{\psi}_{\omega}(m)\overline{\tilde{\sigma}_l(m)}\right)^2\nonumber\\
&=&C\sum^{\bf r}_{l=1}\frac{\theta^2_l\pi_{l}}{N^{\alpha_l}}\left\{\frac{2^{-2\nu_lj}}{\left[\sum^{\bf r}_{l=1}\theta_l2^{-2\nu_lj}\right]^2} \right\}\left(\sigma_{l1} 2^{-j/2}\right)^2
\eeqn
Hence, choosing $\theta_l$ according to \fr{theta1-lm-unif} and combining \fr{delta1-iid}, \fr{delt2} and \fr{delt'22} completes the proof of \fr{var-bet-nonunif-lm}. $\Box$\\
In addition, notice that $z^{(l)}_i=\sum^N_{i=1}\varepsilon^{(l)}_{i, i-1}$ is a zero mean Gaussian. Then, by the Gaussian property $E[z_1^4]=3\left[E|z_1|^2\right]^2$, the fourth moment of $\Delta_{22}$ satisfies
\beqn \label{delt'22-4}
\EE\left[\Delta_{22}\right]^4&\leq&\frac{1}{N^4}\sum^{\bf r}_{l=1}\EE^4\left[\frac{U^{(l)}_{\omega}(t_l)\sigma_l(t_l)}{h_l(t_l)}\right]3\left[\pi_lN^{2-\alpha_l}\right]^2
+\frac{1}{N^4}\left\{\sum^{\bf r}_{l=1}\EE^2\left[\frac{U^{(l)}_{\omega}(t_l)\sigma_l(t_l)}{h_l(t_l)}\right]\pi_lN^{2-\alpha_l}\right\}^2\nonumber\\
&=&C\sum^{\bf r}_{l=1}\frac{\theta^4_l\pi^2_{l}\sigma^4_{l1}}{N^{2\alpha_l}}\left\{\frac{2^{-4\nu_lj}}{\left[\sum^{\bf r}_{l=1}\theta_l2^{-2\nu_lj}\right]^4} \right\}2^{-2j}
+C\left\{\sum^{\bf r}_{l=1}\frac{\theta^2_l\pi_{l}\sigma^2_{l1}}{N^{\alpha_l}}\left\{\frac{2^{-2\nu_lj}}{\left[\sum^{\bf r}_{l=1}\theta_l2^{-2\nu_lj}\right]^2} \right\}2^{-j}\right\}^2\nonumber\\
&\leq& C\left\{\sum^{\bf r}_{l=1}\frac{\theta^2_l\pi_{l}\sigma^2_{l1}}{N^{\alpha_l}}\left\{\frac{2^{-2\nu_lj}}{\left[\sum^{\bf r}_{l=1}\theta_l2^{-2\nu_lj}\right]^2} \right\}2^{-j}\right\}^2.
\eeqn
To prove \fr{var2-bet-nonunif-lm} we combine \fr{D14},  \fr{delt24} and \fr{delt'22-4}, and choose $\theta_l$ according to  \fr{theta1-lm-unif}. $\Box$\\
{\bf Proof of Lemma \ref{lem:Lardev-nonunif-iid}.} Recall \fr{alef-iid}, and notice that
\beqn \label{p-decom}
\Pr \left(| \widehat{\beta_{\omega}}-\beta_{\omega} |> 1/2\lambda({j})\right)&\leq&\Pr\left(|\Delta_1|>1/4\lambda({j})\right)+\Pr\left(|\Delta_2|>1/4\lambda({j})\right)\nonumber\\
&=&P_1+P_2.
\eeqn
For $P_1$, we apply the $i.i.d.$ version of Bernstein inequality described below. 
   \begin{lemma}(Bernstein Inequality: i.i.d. Version). \label{lem:bernineqiid}
Let $\mathcal{Y}_i$, $i=1, 2, \cdots, N$,  be $i.i.d.$ random variables with mean zero and finite variance $\sigma^2$, and let $\|\mathcal{Y}_i\| \leq \|\mathcal{Y}\|_{\infty} < \infty$. Then,  
\be \label{prob-biid}
\Pr\left(\left|N^{-1}\sum^N_{i=1}\mathcal{Y}_i\right|>z\right) \leq 2 \exp\left\{-\frac{Nz^2}{2(\sigma^2+\|\mathcal{Y}\|_{\infty}z/3)}\right\}.
\ee
\end{lemma}
For $\Delta_{1}$, notice that the quantities 
\beqn
\mathcal{Y}_{i}=\sum^{\bf r}_{l=1}\left[\frac{U^{(l)}_{\omega}(t_{li})\sigma_l(t_{li})}{h_l(t_{li})}-\EE\left[\frac{U^{(l)}_{\omega}(t_l)\sigma_l(t_l)}{h_l(t_l)}\right]\right],
\eeqn
are zero-mean $i.i.d.$ random variables with variance satisfying
\beqn
\EE\left[\mathcal{Y}_{i}\right]^2\leq  \|f\|^2_{\infty}\max_l\|g_l\|^2_{\infty}}{\frac{\max_l\left\{ K_{l2}\left(\frac{2\pi}{3}\right)^{-2\nu_l}\right\}\sum^{\bf r}_{l=1}m^{-1}_{l1}\theta^2_l2^{-2\nu_lj}}{\left[\min_l \left\{K_{l1}\left(\frac{8\pi}{3}\right)^{-2\nu_l}\right\}\sum^{\bf r}_{l=1}\theta_l2^{-2\nu_lj}\right]^2},
\eeqn
\beqn
\|\mathcal{Y}_{i}\|_{\infty}\leq 2^{j/2}\|f\|_{\infty}\frac{\max_l\left\{ \|g_l\|_{\infty}\sqrt{K_{l2}}\left(\frac{2\pi}{3}\right)^{-\nu_l}\right\}\sum^{\bf r}_{l=1}m^{-1}_{l1}\theta_l2^{-\nu_lj}}{\left[\min_l \left\{K_{l1}\left(\frac{8\pi}{3}\right)^{-2\nu_l}\right\}\sum^{\bf r}_{l=1}\theta_l2^{-2\nu_lj}\right]}.
\eeqn 
Therefore, by \fr{prob-biid}, as  $N\rightarrow \infty$, and using the fact that $2^j\leq\left[\frac{N}{\ln(N)}\right]^{\frac{1}{2}}$, under $i.i.d.$ errors and long-memory errors with $\Delta=0$, with the choice of $\theta_l$ given by \fr{theta1-iid}, one has
\beqn \label{p1-iid}
P_1&\leq& \Pr\left(\left|\frac{1}{N}\sum^N_{i=1}\mathcal{Y}_{i}\right|>\frac{1}{4}\lambda\right)\nonumber\\
&\leq& 2\exp\left\{-\frac{\frac{N\rho^2\ln(N)}{N\sum^{\bf r}_{l=1}2^{-\nu_lj}}}{32\left(\|f\|^2_{\infty}\frac{\max_l\left\{ \|g_l\|^2_{\infty}{K_{l2}}\left(\frac{2\pi}{3}\right)^{-2\nu_l}\right\}\sum^{\bf r}_{l=1}m_{l1}2^{-2\nu_lj}}{\left[\min_l \left\{K_{l1}\left(\frac{8\pi}{3}\right)^{-2\nu_l}\right\}\sum^{\bf r}_{l=1}m_{l1}2^{-2\nu_lj}\right]^2}+\frac{2^{j/2}\|f\|_{\infty}}{12}\frac{\max_l\left\{ \|g_l\|_{\infty}\sqrt{K_{l2}}\left(\frac{2\pi}{3}\right)^{-\nu_l}\right\}\sum^{\bf r}_{l=1}2^{-\nu_lj}}{\left[\min_l \left\{K_{l1}\left(\frac{8\pi}{3}\right)^{-2\nu_l}\right\}\sum^{\bf r}_{l=1}m_{l1}2^{-2\nu_lj}\right]}\lambda\right)}\right\}\nonumber\\
&\leq& 2\exp\left\{-\frac{{\rho^2\ln(N)}}{32\left(\|f\|^2_{\infty}\frac{\max_l\left\{ \|g_l\|^2_{\infty}{K_{l2}}\left(\frac{2\pi}{3}\right)^{-2\nu_l}\right\}}{\left[\min_l \left\{K_{l1}\left(\frac{8\pi}{3}\right)^{-2\nu_l}\right\}\sqrt{\min_l\{m_{l1}\}}\right]^2}+\frac{2^{j/2}\rho\|f\|_{\infty}}{12}\frac{\max_l\left\{ \|g_l\|_{\infty}\sqrt{K_{l2}}\left(\frac{2\pi}{3}\right)^{-\nu_l}\right\}}{\left[\min_l \left\{K_{l1}\left(\frac{8\pi}{3}\right)^{-2\nu_l}\right\}\min_l\{m_{l1}\}\right]}\sqrt{\frac{\ln(N)}{N}}\right)}\right\}\nonumber\\
&\leq& 2\exp\left\{-\frac{{\rho^2\min_l \left\{K^2_{l1}\left(\frac{8\pi}{3}\right)^{-4\nu_l}\right\}\min_l\{m_{l1}\}\ln(N)}}{64\|f\|^2_{\infty}\max_l\left\{ \|g_l\|^2_{\infty}{K_{l2}}\left(\frac{2\pi}{3}\right)^{-2\nu_l}\right\}}\right\}.
\eeqn
The same way one can show that result \fr{p1-iid} applies to the long-memory errors case under property $\Delta\neq0$, with $\theta_l$ chosen as \fr{theta1-lm-unif}. To evaluate $P_2$ under $i.i.d.$ errors, notice that conditional on ${\bf t}=(t_{11}, t_{12}, \cdots, t_{1N}, t_{21}, t_{22}, \cdots, t_{{\bf r}N})$, $\Delta_2$ is a Gaussian random variable with mean 0 and variance
\be
\Var\left(\Delta_2\mid {\bf t}\right)=\frac{1}{N^2}\sum^{\bf r}_{l=1}\sum^N_{i=1}{\left(\frac{U^{(l)}_{\omega}(t_{li})}{h_l(t_{li})}\right)^2\sigma_l^2(t_{li})}.
\ee
Now, consider the set $\Omega_{\mu}$ defined as
\beqns
\Omega_{\mu}=\left\{{\bf t}: \left|\frac{1}{N}\sum^N_{i=1}\sum^{\bf r}_{l=1}{\left|\frac{U^{(l)}_{\omega}(t_{li})}{h_{l}(t_{li})}\right|^2}-\EE\left[\sum^{\bf r}_{l=1}{\left|\frac{U^{(l)}_{\omega}(t)}{h_l(t)}\right|^2}\right]\right|>\mu\frac{\sum^{\bf r}_{l=1}m^{-1}_{l1}\theta^2_l2^{-2\nu_lj}}{\left[\sum^{\bf r}_{l=1}\theta_l2^{-2\nu_lj}\right]^2}\right\}.
\eeqns
In addition, notice that for ${\bf t}\in \Omega^c_{\mu}$, we have 
\beqn \label{omeg-mu-b}
\frac{1}{N}\sum^N_{i=1}\sum^{\bf r}_{l=1}{\left|\frac{U^{(l)}_{\omega}(t_{li})}{h_l(t_{li})}\right|^2}&\leq& \left(\EE\left[\sum^{\bf r}_{l=1}{\left|\frac{U^{(l)}_{\omega}(t)}{h_l(t)}\right|^2}\right]+\left|\frac{1}{N}\sum^N_{i=1}\sum^{\bf r}_{l=1}{\left|\frac{U^{(l)}_{\omega}(t_{li})}{h_l(t_{li})}\right|^2}-\EE\left[\sum^{\bf r}_{l=1}{\left|\frac{U^{(l)}_{\omega}(t)}{h_l(t)}\right|^2}\right]\right|\right)\nonumber\\
&\leq& \left(\EE\left[\sum^{\bf r}_{l=1}{\left|\frac{U^{(l)}_{\omega}(t)}{h_l(t)}\right|^2}\right]+\mu\frac{\sum^{\bf r}_{l=1}m^{-1}_{l1}\theta^2_l2^{-2\nu_lj}}{\left[\sum^{\bf r}_{l=1}\theta_l2^{-2\nu_lj}\right]^2}\right). 
\eeqn
Then, $P_{2}$ can be partitioned as $P_{2}\leq P_{21}+P_{22}$,  where
\beqn
P_{21}&=&\EE\left[\Pr\left(\left|\Delta_{2}\right|>1/4\lambda(j)\mid {\bf t}\right)\II\left(\Omega^c_{\mu}({\bf t}\right)\right].\\
P_{22}&=&\EE\left[\II\left(\Omega_{\mu}({\bf t})\right)\right].
\eeqn
Hence, with $\theta_l$ chosen as \fr{theta1-iid}, by Gaussian tail inequality and result \fr{omeg-mu-b}, one has
\beqn \label{p21-iid}
P_{21}&\leq& 2\exp\left\{-\frac{\rho^2\ln(N)\left[\sum^{\bf r}_{l=1}m_{l1}2^{-2\nu_l j}\right]^{-1}}{32\left(\frac{\max_l\left\{\sigma^2_{l2}K_{l2}\left(\frac{2\pi}{3}\right)^{-2\nu_l}\right\}\sum^{\bf r}_{l=1}m_{l1}2^{-2\nu_lj}}{\left[\min\left\{K_{l1}\left(\frac{8\pi}{3}\right)^{-2\nu_l}\right\}\sum^{\bf r}_{l=1}m_{l1}2^{-2\nu_lj}\right]^2}+\mu\frac{\sum^{\bf r}_{l=1}m_{l1}2^{-2\nu_lj}}{\left[\sum^{\bf r}_{l=1}m_{l1}2^{-2\nu_lj}\right]^2}\right)}\right\}\nonumber\\
&=& N^{-\frac{\min\left\{K^2_{l1}\left(\frac{8\pi}{3}\right)^{-4\nu_l}\right\}\rho^2}{32\max_l\left\{\sigma^2_{l2}K_{l2}\left(\frac{2\pi}{3}\right)^{-2\nu_l}\right\}+32\mu \min_l\left\{K^2_{l1}\left(\frac{8\pi}{3}\right)^{-4\nu_l}\right\}}}.
\eeqn
For $P_{22}$, notice that {\bf Lemma \ref{lem:bernineqiid}} can be applied with 
\beqn
\mathcal{Y}_{i}=\left[\sum^{\bf r}_{l=1}{\left|\frac{U^{(l)}_{\omega}(t_{li})}{h_l(t_{li})}\right|^2}-\EE\left[\sum^{\bf r}_{l=1}{\left|\frac{U^{(l)}_{\omega}((t)}{h_l(t)}\right|^2}\right]\right],
\eeqn
\beqn \label{var-y}
\EE\left[\mathcal{Y}_{i}\right]^2\leq 8\frac{\max_l\left\{ K^2_{l2}\left(\frac{2\pi}{3}\right)^{-4\nu_l}\right\}}{\min_l \left\{K^4_{l1}\left(\frac{8\pi}{3}\right)^{-8\nu_l}\right\}}\left(\frac{\left[\sum^{\bf r}_{l=1}\theta^2_l2^{-2\nu_lj}\right]^2}{\left[\sum^{\bf r}_{l=1}\theta_l2^{-2\nu_lj}\right]^4}\right),
\eeqn
\beqn \label{max-y}
\|\mathcal{Y}_{i}\|_{\infty}\leq \frac{\max_l\left\{ K_{l2}\left(\frac{2\pi}{3}\right)^{-2\nu_l}\right\}}{\min_l \left\{K^2_{l1}\left(\frac{8\pi}{3}\right)^{-4\nu_l}\right\}}\left(\frac{\sum^{\bf r}_{l=1}\theta^2_l2^{-2\nu_lj}}{\left[\sum^{\bf r}_{l=1}\theta_l2^{-2\nu_lj}\right]^2}\right).
\eeqn 
Hence, using $\theta_l=m_{l1}$, by Bernstein inequality, 
\beqn \label{p22-iid}
P_{22}&\leq& 2\exp\left\{-\frac{N\mu^2\frac{1}{\left[\sum^{\bf r}_{l=1}m_{l1}2^{-2\nu_lj}\right]^2}}{2\left(8\frac{\max_l\left\{ K^2_{l2}\left(\frac{2\pi}{3}\right)^{-4\nu_l}\right\}}{\min_l \left\{K^4_{l1}\left(\frac{8\pi}{3}\right)^{-8\nu_l}\right\}}\left(\frac{\max_l\{m_{l1}\}}{\left[\sum^{\bf r}_{l=1}m_{l1}2^{-2\nu_lj}\right]^2}\right)+\frac{\mu}{3}\frac{\max_l\left\{ K_{l2}\left(\frac{2\pi}{3}\right)^{-2\nu_l}\right\}}{\min_l \left\{K^2_{l1}\left(\frac{8\pi}{3}\right)^{-4\nu_l}\right\}}\left(\frac{\sqrt{\max_l\{m_{l1}\}}}{\sum^{\bf r}_{l=1}m_{l1}2^{-2\nu_lj}}\right)^2\right)}\right\}\nonumber\\
&=& 2\exp\left\{-\frac{N\mu^2}{2\left(8\frac{\max_l\left\{ K^2_{l2}\left(\frac{2\pi}{3}\right)^{-4\nu_l}\right\}\max_l\{m_{l1}\}}{\min_l \left\{K^4_{l1}\left(\frac{8\pi}{3}\right)^{-8\nu_l}\right\}}+\frac{\mu}{3} \frac{\max_l\left\{ K_{l2}\left(\frac{2\pi}{3}\right)^{-2\nu_l}\right\}\max_l\{m_{l1}\}}{\min_l \left\{K^2_{l1}\left(\frac{8\pi}{3}\right)^{-4\nu_l}\right\}}\right)}\right\}\nonumber\\
&=&o\left( N^{-2\tau}\right). 
\eeqn
Hence, combining \fr{p1-iid}, \fr{p21-iid}, and \fr{p22-iid} and choosing 
\beqn\label{rho-bound}
\rho^2\geq 18\tau \max\left\{\frac{\|f\|^2_{\infty}\max_l\left\{ K^2_{l2}\left(\frac{2\pi}{3}\right)^{-4\nu_l}\right\}}{\min_l \left\{K^2_{l1}m_{l1}\left(\frac{8\pi}{3}\right)^{-4\nu_l}\right\}}, \left[\frac{\max_l\left\{ K^2_{l2}\left(\frac{2\pi}{3}\right)^{-4\nu_l}\right\}}{\min_l \left\{K^2_{l1}\left(\frac{8\pi}{3}\right)^{-4\nu_l}\right\}}\vee \mu\right]\right\},
\eeqn
completes the proof of \fr{Largdev-iid} under $i.i.d.$ errors.  $\Box$\\
 In order to prove \fr{Largdev-iid} under long-memory errors, we make use of martingale version of Bernstein inequality. 
    \begin{lemma}(Bernstein Inequality: Martingale Version). \label{lem:bernineqmar}
Let $\left(d_i, \mathcal{F}_i\right)$, $i=1, 2, \cdots, N$,  be a martingale difference sequence.  Then,  
\be \label{prob-bmar}
\Pr\left(\left|\sum^N_{i=1}d_i\right|>z, \sum^N_{i=1}d_i^2\II(|d_i|>x)+\sum^N_{i=1}\EE\left[d_i^2\mid \mathcal{F}_{i-1}\right]\leq L\right) \leq 2 \exp\left\{-\frac{z^2}{2(L+xz/3)}\right\}.
\ee
\end{lemma}
Notice that under long-memory errors, $P_2$  in \fr{p-decom} can be decomposed according to martingale and long-memory term as
\beqn
P_2&\leq&\Pr\left(\left|\Delta_{21}\right|>1/8\lambda(j)\right)\nonumber\\
&+& \Pr\left(\left|\Delta_{22}\right|>1/8\lambda(j)\right)=P_{21}+P_{22}. 
\eeqn
In case $\Delta=0$, $\Delta_{22}$ defined in \fr{delta2} is such that $\Delta_{22}=0$, and therefore, $P_{22}=0$.  However, in case $\Delta\neq0$, then  $\Delta_{22}$ is a zero-mean Gaussian random variable with variance 
\beqn \label{var-delt22}
\Var\left[\Delta_{22}\right]&\leq&\sum^{\bf r}_{l=1}\frac{\theta^2_l\pi_{l}2^{-j}}{N^{\alpha_l}}\left\{\frac{\sigma^2_{l1}K_{l2}\left(\frac{2\pi}{3}\right)^{-2\nu_l}2^{-2\nu_lj}}{\left[\sum^{\bf r}_{l=1}\theta_lK_{l1}\left(\frac{8\pi}{3}\right)^{-2\nu_l}2^{-2\nu_lj}\right]^2} \right\}\nonumber\\
&\leq &\left\{\frac{\sum^{\bf r}_{l=1}\theta_l\sigma^2_{l1}K_{l2}\left(\frac{2\pi}{3}\right)^{-2\nu_l}2^{-2\nu_lj}}{\left[\sum^{\bf r}_{l=1}\theta_lK_{l1}\left(\frac{8\pi}{3}\right)^{-2\nu_l}2^{-2\nu_lj}\right]^2} \right\}\nonumber\\
&\leq&\frac{\max_l\left\{ \sigma^2_{l1}K_{l2}\left(\frac{2\pi}{3}\right)^{-2\nu_l}\right\}}{\min_l \left\{K^2_{l1}\left(\frac{8\pi}{3}\right)^{-4\nu_l}\right\}}\frac{1}{\sum^{\bf r}_{l=1}\theta_l2^{-2\nu_lj}}.
\eeqn
Hence, by the Gaussian tail probability inequality, thresholds \fr{thresh-lm1-nonunif} and result \fr{var-delt22}, $P_{22}$ satisfies
\beqn\label{p22-lm}
P_{22}&\leq& 2 \exp\left\{-\frac{\rho^2\min_l \left\{K^2_{l1}\left(\frac{8\pi}{3}\right)^{-4\nu_l}\right\}\ln^2(N)}{128\max_l\left\{ \sigma^2_{l1}K_{l2}\left(\frac{2\pi}{3}\right)^{-2\nu_l}\right\}}\right\}\nonumber\\
&=&o\left(N^{-2\tau}\right). 
\eeqn
For $P_{21}$, denote 
\beqns
\Lambda^{(l)}(N, z_l)=\sum^N_{i=1}(d^{(l)}_i)^2\II(|d^{(l)}_i|>z_l)+\sum^N_{i=1}\EE\left[(d^{(l)}_i)^2\mid \mathcal{F}^{(l)}_{i-1}\right],
\eeqns
 and notice that
\beqn \label{p2}
P_{21}&=&\Pr\left(|\Delta_{21}|> \frac{1}{8} \lambda \right)\leq \Pr\left(\frac{1}{N}\sum^{\bf r}_{l=1}\left|\sum^N_{i=1}d^{(l)}_i\right|>\frac{1}{8}\lambda \right)\nonumber\\
&\leq&\sum^{\bf r}_{l=1}\left[ \Pr\left( \frac{1}{N}\left|\sum^N_{i=1}d^{(l)}_i\right|>\frac{1}{8{\bf r}}\lambda, \Lambda^{(l)}(N, z_l)\leq L^{(l)} \right)+\Pr\left(\Lambda^{(l)}(N, z_l)>L^{(l)}\right)\right].
\eeqn
Thus, under property $\Delta= 0$, by result \fr{prob-bmar}  and  as $N \rightarrow \infty$, the first summation term of \fr{p2} has components that satisfy
\beqn\label{p21}
\Pr\left( \frac{1}{N}\left|\sum^N_{i=1}d^{(l)}_i\right|>\frac{1}{8{\bf r}}\lambda, \Lambda^{(l)}(N, z_l)\leq L^{(l)} \right)\leq 2\exp\left\{\frac{-N^2\ln^2(N)\rho^2}{128{\bf r}^2\sum^{\bf r}_{l=1}\theta_l2^{-2\nu_lj}\left(L^{(l)}+\frac{z_lN\rho}{24{\bf r}}\frac{\ln(N)}{\sqrt{\sum^{\bf r}_{l=1}\theta_l2^{-2\nu_lj}}}\right)}\right\}. \nonumber\\
\eeqn
As for the second term of \fr{p2}, at each $l$ value, it can be decomposed into two terms as follows
\beqn\label{p31}
\Pr\left(\Lambda^{(l)}(N, z_l)>L^{(l)}\right)&\leq& \left[\Pr\left(\sum^N_{i=1}\left[d^{(l)}_i\right]^2\II(|d^{(l)}_i|>z_l)> \frac{L^{(l)}}{2}\right)+\Pr\left(\sum^N_{i=1}\EE\left[\left[d^{(l)}_i\right]^2\mid \mathcal{F}^{(l)}_{i-1}\right]>\frac{L^{(l)}}{2}\right)\right]\nonumber\\
&=&P_{2121}(l)+P_{2122}(l).
\eeqn
For $P_{2122}(l)$, notice that by choosing 
\be\label{Ll}
L^{(l)}=\frac{\sigma^2_{l1}K_{l2}m^{-1}_{l1}\theta^2_l2^{-2\nu_l j}\left(\frac{2\pi}{3}\right)^{-2\nu_l}}{\left[\sum^{\bf r}_{l=1}K_{l1}\theta_l\left(\frac{8\pi}{3}\right)^{-2\nu_l}2^{-2\nu_l j}\right]^2}2N\left[D_l\ln(N)+\EE\left[\eta^{(l)}_1\right]^2\right],
\ee
 and keeping in mind that $\varepsilon^{(l)}_{i, i-1}$ are Gaussian, one has 
\beqn \label{p322}
P_{2122}(l)= \Pr\left(\sum^N_{i=1}\EE\left[\left[d^{(l)}_i\right]^2\mid \mathcal{F}^{(l)}_{i-1}\right]>L/2\right)&\leq& \Pr\left(\sum^N_{i=1}\left[\varepsilon^{(l)}_{i, i-1}\right]^2> D_lN\ln(N)\right).\nonumber\\
&\leq & 2N\exp\left\{-\frac{D_l\ln(N)}{2\EE[\varepsilon^{(l)}_{1,0}]^2}\right\}
\eeqn
As for $P_{2121}(l)$, notice that $2^{j}\leq\frac{{N}^{\frac{1}{2\nu_o+1}}}{\ln(N)}$. Since $\varepsilon^{(l)}_{i, i-1}$ and $\eta^{(l)}_i$ are Gaussians and independent from each other, taking 
\beqn\label{zl}
z^2_l=\frac{c_{lo}{N}^{\frac{1}{2\nu_o+1}}\theta_l2^{-2\nu_l j}}{\left[\sum^{\bf r}_{l=1}\theta_l2^{-2\nu_l j}\right]^2},
\eeqn
as $N \rightarrow \infty$, one has
\beqn \label{p321}
P_{2121}(l)&=& \Pr\left(\sum^N_{i=1}\left[d_i\right]^2\II(|d^{(l)}_i|>z_l)> \frac{L^{(l)}}{2}\right)\nonumber\\
&\leq& N \Pr\left(\left[d^{(l)}_1\right]^2\II(|d^{(l)}_1|>z_l)>\frac{\sigma^2_{l1}K_{l2}m^{-1}_{l1}\theta^2_l\left(\frac{2\pi}{3}\right)^{-2\nu_l}2^{-2\nu_l j}}{\left[\sum^{\bf r}_{l=1}\theta_lK_{l1}\left(\frac{8\pi}{3}\right)^{-2\nu_l}2^{-2\nu_l j}\right]^2}D_l\ln(N)\right)\nonumber\\
&\leq & N\Pr\left(\left[d^{(l)}_1\right]^2>\max\left\{\frac{\sigma^2_{l1}K_{l2}m^{-1}_{l1}\theta^2_l\left(\frac{2\pi}{3}\right)^{-2\nu_l}2^{-2\nu_l j}}{\left[\sum^{\bf r}_{l=1}\theta_lK_{l1}\left(\frac{8\pi}{3}\right)^{-2\nu_l}2^{-2\nu_l j}\right]^2}D_l\ln(N), z_l^2\right\}\right)\nonumber\\
&\leq & N\left[\Pr\left(\frac{\sigma^2_{l1}2^jK_{l2}m^{-1}_{l1}\theta^2_l\left(\frac{2\pi}{3}\right)^{-2\nu_l}2^{-2\nu_l j}}{\left[\sum^{\bf r}_{l=1}\theta_lK_{l1}\left(\frac{8\pi}{3}\right)^{-2\nu_l}2^{-2\nu_l j}\right]^2}2\left[|\varepsilon^{(l)}_{1, 0}|^2+|\eta^{(l)}_{1}|^2\right]>{z_l^2}\right)\right]\nonumber\\
&\leq& 2N\exp\left\{\frac{-\left[\sum^{\bf r}_{l=1}\theta_lK_{l1}\left(\frac{8\pi}{3}\right)^{-2\nu_l}2^{-2\nu_l j}\right]^2z_l^2}{8 \sigma^2_{l1}K_{l2}m^{-1}_{l1}\theta^2_l\left(\frac{2\pi}{3}\right)^{-2\nu_l}2^{-2\nu_l j}2^{j}\left(\EE\left[\varepsilon^{(l)}_{1,0}\right]^2\wedge\EE\left[\eta^{(l)}_{1}\right]^2\right)}\right\}\nonumber\\
& \leq & 2N\exp\left\{-\frac{c_{lo}\min_l\{K_{l1}^2\left(\frac{3}{8\pi}\right)^{4\nu_l}\}\ln(N)}{8\sigma^2_{l1}K_{l2}\left(\frac{3}{2\pi}\right)^{2\nu_l}\left(\EE\left[\varepsilon^{(l)}_{1,0}\right]^2\wedge\EE\left[\eta^{(l)}_{1}\right]^2\right)}\right\}.
\eeqn
Notice that with the choices of $L_l$ and $z_l$ given by \fr{Ll} and \fr{zl}, \fr{p21} will be of order $o\left( \left[\frac{1}{N}\right]^{2\tau}\right)$. To complete the proof we choose $D_l$ and $c_{lo}$ such that
\beqns
D_l\geq 2(1+2\tau)\EE\left[\varepsilon^{(l)}_{1,0}\right]^2, \ \ c_{lo}\geq 8(1+2\tau)\frac{\sigma_{l1}^2K_{l2}\left(\frac{2\pi}{3}\right)^{-2\nu_l}}{\min_l\{K_{l1}^2\left(\frac{8\pi}{3}\right)^{-2\nu_l}\}}\max\left\{\EE[\varepsilon^{(l)}_{1,0}]^2, \EE[\eta^{(l)}_{1}]^2\right\}. \Box
\eeqns
In a similar fashion, one can prove results \fr{p21}, \fr{p322} and \fr{p321} under property $\Delta\neq 0$, so we skip it. $\Box$\\
{\bf Proof of Theorem \ref{th:upperbds-2-nonunif-iid}}. For the $i.i.d.$ errors and long-memory errors with $\Delta=0$ scenarios, denote
\be\label{bet-iid}
2^{J_{o}}=[\chi(A)]^{-\frac{d}{2s'}},
\ee 
and
\be \label{chijj-iid}
\chi(A)=\frac{\ln(N)}{A^2N},
\ee
 where,
\be
d_o=\min\left\{\frac{2s'}{2s'+2\nu_o}, \frac{2s}{2s+2\nu_o+1}\right\}.
\ee 
Note that with the choices of $J$ and $\lambda(j)$ given by \fr{finresol}  and \fr{thresh-lm1-nonunif-iid}, respectively, the estimation error can be decomposed as $\mathbb{E} \| \widehat{f}_N-f \|^2
\leq\mathbb{E}_1 +\mathbb{E}_2 +\mathbb{E}_3+\mathbb{E}_4+\mathbb{E}_5+\mathbb{E}_6$, where
 \beqn  
\mathbb{E}_1&=&  \sum^{J-1}_{j=m_{o}}\sum^{2^{j}-1}_{k=0}\EE \left[\left| \widehat{ \beta_{\omega}}-\beta_{\omega}\right|^2 \II \left(   \left| \widehat{ \beta_{\omega}}-\beta_{\omega}  \right| > \frac{1}{2} \lambda\right)  \right], \label{r21}\\
\mathbb{E}_2&=&  \sum^{J-1}_{j=m_{o}}\sum^{2^{j}-1}_{k=0}\EE \left[\left| \widehat{ \beta_{\omega}}-\beta_{\omega}\right|^2 \II \left(  \left|  \beta_{\omega} \right| >  \frac{1}{2}\lambda \right)  \right],\label{r22}\\
\mathbb{E}_3&=&  \sum^{J-1}_{j=m_{o}}\sum^{2^{j}-1}_{k=0}\beta_{\omega}^2 \Pr \left( \left| \widehat{ \beta_{\omega}}-\beta_{\omega}  \right| > \frac{1}{2} \lambda \right),\label{r31}\\
\mathbb{E}_4&=& \sum^{J-1}_{j=m_{o}}\sum^{2^{j}-1}_{k=0}\beta_{\omega}^2\II \left(  \left|  \beta_{\omega} \right| <  \frac{3}{2}\lambda \right),\label{r32}\\
\mathbb{E}_5&=& \sum^{\infty}_{j=J}\sum^{2^{j}-1}_{k=0} \left(\beta^*_{\omega}\right)^2,\label{r33}
\eeqn
\beqn  
\mathbb{E}_6&=&\sum^{2^{m_{o}}-1}_{k=0} \EE \left[\left| \widehat{\beta}_{m_{o}, k} - \beta_{m_{o}, k}\right|^2\right]. \label{r34}
\eeqn
Then, by  \fr{var-bet-nonunif-iid}, as $N\rightarrow \infty$, \fr{r34} becomes,
\beqn \label{e6-iid}
\mathbb{E}_6=O\left(\frac{1}{N}\right). 
\eeqn
Now, by \fr{eq11} and \fr{finresol}, as $N\rightarrow \infty$, \fr{r33} becomes 
\beqn \label{ee14-iid}
\mathbb{E}_5=O\left( \sum^{\infty}_{j=J}A^2 2^{-j2s'} \right)=O\left(A^2 \left[\frac{1}{A^2N}\right]^{\frac{2s'}{2\nu_o+1}}\right)
=O\left(A^2 \left[\chi(A)\right]^{d_o}\right).
\eeqn
Now, combining $\mathbb{E}_1$ and $\mathbb{E}_3$, and applying Cauchy-Schwarz inequality, the moments property of the Gaussian, {\bf Lemma \ref{lem:Var-nonunif}} with the choice $\mu<\frac{\max_l\left\{ K^2_{l2}\left(\frac{2\pi}{3}\right)^{-4\nu_l}\right\}}{\min_l \left\{K^2_{l1}\left(\frac{8\pi}{3}\right)^{-4\nu_l}\right\}}$ in \fr{rho-bound}, \fr{eq11} and \fr{finresol}, as $N\rightarrow \infty$, yields 
\beqn
\mathbb{E}_1 +\mathbb{E}_3&=&O\left(\frac{\left[\sum^{\bf r}_{l=1}2^{-(2\nu_l+1)J}\right]^{-1}}{N}\left[ \frac{1}{N}\right]^{{\tau}}+ A^22^{2-m_{o}s'}\left[ \frac{1}{N}\right]^{{2\tau}} \right)\nonumber\\
&=&O\left(\frac{2^{(2\nu_o+1)J}}{N}\left[ \frac{1}{N}\right]^{{\tau}}+\left[ \frac{1}{N}\right]^{{2\tau}} \right)=  O \left(  \left[ \frac{1}{N}\right]^{{\tau}} \right). \label{r21r31-iid}
\eeqn
Now, combining $\mathbb{E}_2$ and $\mathbb{E}_4$ yields
\beqn  \label{r22r32}
\Delta=\mathbb{E}_2 +\mathbb{E}_4 &=&\sum^{J-1}_{j=m_{o}}\sum^{2^{j}-1}_{k=0}\left[\min\left\{\beta_{\omega}^2, \left[\sum^{\bf r}_{l=1}2^{-2\nu_lj}\right]^{-1}\frac{\ln(N)}{N} \right\} \right]\nonumber\\
 &=&\sum^{J-1}_{j=m_{o}}\sum^{2^{j}-1}_{k=0}\left[\min\left\{\beta_{\omega}^2, 2^{2\nu_o j}\frac{\ln(N)}{N} \right\} \right].
\eeqn
This latter can be partitioned into 
\beqn  
 \Delta_{1}&=&  O \left(\sum^{J-1}_{j=J_{o}}A^2 2^{-2js'} \right)= O\left(A^2 \left[\chi(A)\right]^{d_o}\right), \label{del1-iid}\\
 \Delta_{2}&=&O \left( \sum^{J_{o}-1}_{j=m_{o}}\frac{ 2^{j2\nu_o+j}\ln(N)}{N}\II \left(  2^{j2\nu_o+j}\leq   \left[ \chi(A)\right]^{d-1} \right)\right)=O\left(A^2 \left[ \chi(A)\right]^{d_o}\right), \label{del2-iid}\\
 \Delta_{3}&=&O \left( \sum^{J_{o}-1}_{j=m_{o}}A^{p'}2^{-p'js'}\left[\frac{\ln(N)}{N}2^{2\nu_o j}\right]^{1-\frac{p'}{2}}\II \left(  2^{j2\nu_o+j}>  \left[ \chi(A)\right]^{d_o-1} \right)\right). \label{del3-iid}
\eeqn
$\Delta_3$ will be studied under different cases as follows:\\
\noindent Case 1: $s>(\frac{1}{p}-\frac{1}{2})(2\nu_o+1)$. Then, $d_o=\frac{2s}{2s+2\nu_o+1}$, and 
for $\Delta_{3}$ in \fr{del3-iid}, one has 
\beqn  
  \Delta_{3}&=&O \left( A^2 \left[ \chi(A)\right]^{1-\frac{p'}{2}}\sum^{J_{o}-1}_{j=m_{o}}2^{-j(p's'-(1-\frac{p'}{2})2\nu_o)}\II \left(  2^{j(2\nu_o+1)}>  \left[ \chi(A)\right]^{d_o-1} \right)\right).\nonumber\\
  &=& O \left( A^2 \left[ \chi(A)\right]^{1-\frac{p'}{2}}\sum^{J_{o}-1}_{j=m_{o}}2^{-jp'(s-(\frac{1}{p'}-\frac{1}{2})(2\nu_o+1))}\II \left(  2^{j(2\nu_o+1)}>  \left[ \chi(A)\right]^{d_o-1} \right) \right)\nonumber\\
  &=& O \left( A^2 \left[ \chi(A)\right]^{\frac{2s}{2s+2\nu_o+1}}\right).  \label{del'3-iid}
\eeqn
\noindent Case 2: $s\leq(\frac{1}{p}-\frac{1}{2})(2\nu_o+1)$. Then, $d_o=\frac{2s'_1}{2s'+2\nu_o}$, and 
for $\Delta_{3}$ in \fr{del3-iid}, one has 
\beqn  
  \Delta_{3}&=&O \left( A^2 \left[ \chi(A)\right]^{1-\frac{p'}{2}}\sum^{J_{o}-1}_{j=m_{o}}2^{-j(p's'-(1-\frac{p'}{2})2\nu_o)}\II \left(  2^{j(2\nu_o+1)}>  \left[ \chi(A)\right]^{d_o-1} \right)\right).\nonumber\\
  &=& O \left( A^2 \left[ \chi(A)\right]^{1-\frac{p'}{2}}\sum^{J_{o}-1}_{j=m_{o}}2^{jp'((\frac{1}{p}-\frac{1}{2})(2\nu_o+1)-s)} \right).\nonumber\\
  &=& O \left( A^2 \left[ \chi(A)\right]^{\frac{2s'}{2s'+2\nu_o}}\left[\ln(N)\right]^{\II\left(s=(\frac{1}{p}-\frac{1}{2})2\nu_o\right)}\right).  \label{del''3-iid}
\eeqn
For the long-memory errors scenario with the property $\Delta=0$, we can add the long-memory phase,$\alpha_o<d_o$, under which 
\be\label{del-lm}
  \Delta_{3}=O \left( A^2  \left[\frac{\ln(N)}{{A}^2N^{\alpha_o}}\right]\right). 
\ee 
Hence, combining \fr{e6-iid}, \fr{ee14-iid}, \fr{r21r31-iid}, \fr{del1-iid}, \fr{del2-iid}, \fr{del'3-iid}, \fr{del''3-iid} completes the proof of \fr{upperbds-2-iid} and combining combining \fr{e6-iid}, \fr{ee14-iid}, \fr{r21r31-iid}, \fr{del1-iid}, \fr{del2-iid}, \fr{del'3-iid}, \fr{del''3-iid}, and \fr{del-lm} completes the proof of \fr{upperbds-2-detl0}. $\Box$\\
For the long memory errors scenario under property $\Delta\neq0$, denote
\be\label{bet}
2^{J_{o}}=[\chi(A)]^{-\frac{d_*}{2s'}},\ \ 2^{J^*_{o}}=[\chi_{\alpha^*}(A)]^{-\frac{d^*}{2s'}}, \ \ 2^{J_1}=N^{1-\alpha^*}.
\ee 
and
\be \label{chijj}
\chi(A)=\frac{\ln^2(N)}{A^2N}, \ \ \chi_{\alpha^*}(A)=\frac{\ln(N)}{A^2N^{\alpha^*}},
\ee
 where,
\be
d_*=\min\left\{\frac{2s'}{2s'+2\nu^*}, \frac{2s}{2s+2\nu^*+1}\right\},\ \ \ d^*=\min\left\{\frac{2s}{2s+2\nu^*}, \frac{2s'}{2s'+2\nu^*-1}\right\},
\ee 
 and note that with the choices of $J$ and $\lambda(j)$ given by \fr{finresol-lm*}  and \fr{thresh-lm1-nonunif}, respectively, the estimation error can be decomposed as $\mathbb{E} \| \widehat{f}_N-f \|^2
\leq\mathbb{E}_1 +\mathbb{E}_2 +\mathbb{E}_3+\mathbb{E}_4+\mathbb{E}_5+\mathbb{E}_6$, according to \fr{r21}, \fr{r22}, \fr{r31}, \fr{r32}, \fr{r33} and \fr{r34}, respectively. Then, by  \fr{var-bet-nonunif-lm}, as $N\rightarrow \infty$, \fr{r34} becomes,
\beqn \label{e6}
\mathbb{E}_6=O\left(\frac{1}{N}+\frac{1}{N^{\alpha^*}} \right). 
\eeqn
Now, by \fr{eq11} and \fr{finresol-lm*}, as $N\rightarrow \infty$, \fr{r33} becomes 
\beqn \label{ee14}
\mathbb{E}_5=O\left( \sum^{\infty}_{j=J}A^2 2^{-j2s'} \right)=O\left(A^2 \left[\frac{1}{A^2N}\right]^{\frac{2s_1'}{2\nu^*+1}}\right)
=O\left(A^2 \left[\chi(A)\right]^{d_*}\vee A^2 \left[\ \chi_{\alpha^*}(A)\right]^{d^*}\right).
\eeqn
Now, combining $\mathbb{E}_1$ and $\mathbb{E}_3$, and applying Cauchy-Schwarz inequality, the moments property of the Gaussian, {\bf Lemma \ref{lem:Var-nonunif}} with the choice $\tau>1$, \fr{eq11} and \fr{finresol-lm*}, as $N\rightarrow \infty$, yields 
\be
\mathbb{E}_1 +\mathbb{E}_3 = O  \left( \frac{2^{(2\nu^*+1) J}}{N}\left[ \frac{1}{N}\right]^{{\tau}}+ A^22^{2-m_{o}s'}\left[ \frac{1}{N}\right]^{{2\tau}} \right)=  O \left(  \left[ \frac{1}{N}\right]^{{\tau}} \right). \label{r21r31}
\ee
To evaluate $\mathbb{E}_2$ and $\mathbb{E}_4$ we need partition according to $\left\{J_1: 2^{J_{1}}=\frac{N}{N^{\alpha^*}}\right\}$. Now, combining $\mathbb{E}_2$ and $\mathbb{E}_4$ yields
\beqn  \label{r22r32}
 \mathbb{E}_2 +\mathbb{E}_4&=&  \sum^{J-1}_{j=m_{o}}\sum^{2^{j}-1}_{k=0} \left[\EE\left| \widehat{ \beta_{\omega}}-\beta_{\omega}\right|^2 \II \left(  \left|  \beta_{\omega} \right| >  \frac{1}{2}\lambda \right)+ \beta^2_{\omega}\II \left(  \left|  \beta_{\omega} \right| <  \frac{3}{2}\lambda \right) \right],\nonumber\\
 &=&\sum^{J_{1}-1}_{j=m_{o}}\sum^{2^{j}-1}_{k=0}\left[\EE\left| \widehat{ \beta_{\omega}}-\beta_{\omega}\right|^2 \II \left(  \left|  \beta_{\omega} \right| >  \frac{1}{2}\lambda \right)+ \beta^2_{\omega}\II \left(  \left|  \beta_{\omega} \right| <  \frac{3}{2}\lambda \right) \right]\nonumber\\
 &+&\sum^{J-1}_{j=J_{1}}\sum^{2^{j}-1}_{k=0} \left[\EE\left| \widehat{ \beta_{\omega}}-\beta_{\omega}\right|^2 \II \left(  \left|  \beta_{\omega} \right| >  \frac{1}{2}\lambda \right)+ \beta^2_{\omega}\II \left(  \left|  \beta_{\omega} \right| <  \frac{3}{2}\lambda \right) \right]\nonumber\\
 &=&\Delta_1+\Delta_2.
\eeqn
Using \fr{thresh-lm1-nonunif} and \fr{finresol-lm*}, as $N\rightarrow \infty$, yields 
\beqn
\Delta_1&=&\sum^{J_{1}-1}_{j=m_{o}}\sum^{2^{j}-1}_{k=0}\left[\min\left\{\beta^2_{\omega}, 2^{2\nu^* j}\frac{2^{-j}\ln(N)}{N^{\alpha^*}} \right\} \right].
 \eeqn
Notice that if $2^{J_1}> 2^{J^*_{o}}$, $\Delta_1$ can be decomposed into the following terms 
\beqn  
 \Delta_{11}&=&  O \left(\sum^{J_{1}-1}_{j=J^*_{o}}A^2 2^{-2js'} \right)= O\left(A^2 \left[\chi_{\alpha^*}(A)\right]^{d^*}\right), \label{del11}\\
 \Delta_{12}&=&O \left( \sum^{J^*_{o}-1}_{j=m_{o}}\frac{ 2^{j2\nu^*}\ln(N)}{N^{\alpha^*}}\II \left(  2^{j2\nu^*}\leq   \left[ \chi_{\alpha^*}(A)\right]^{d^*-1} \right)\right)=O\left(A^2 \left[ \chi_{\alpha^*}(A)\right]^{\frac{2s'}{2s'+2\nu^*}}\right). \label{del12}\\
 \Delta_{13}&=&O \left( \sum^{J^*_{o}-1}_{j=m_{o}}A^{p'}2^{-p'js'}\left[\frac{\ln(N)}{N^{\alpha^*}}2^{2\nu^* j-j}\right]^{1-\frac{p'}{2}}\II \left(  2^{j2\nu^*}>  \left[ \chi_{\alpha^*}(A)\right]^{d^*-1} \right)\right).\nonumber \label{del13}
\eeqn
\noindent Case 1: $s>(\frac{1}{p}-\frac{1}{2})2\nu^*$. Then, $d^*=\frac{2s}{2s+2\nu^*}$, and 
for $\Delta_{13}$ in \fr{del13}, one has 
\beqn  
  \Delta_{13}&=&O \left( A^2 \left[ \chi_{\alpha^*}(A)\right]^{1-\frac{p'}{2}}\sum^{J^*_{o}-1}_{j=m_{o}}2^{-j(p's'-(1-\frac{p'}{2})(2\nu^*-1))}\II \left(  2^{j2\nu^*}>  \left[ \chi_{\alpha^*}(A)\right]^{\frac{2s}{2s+2\nu^*}-1} \right)\right).\nonumber\\
  &=& O \left( A^2 \left[ \chi_{\alpha^*}(A)\right]^{1-\frac{p'}{2}}\sum^{J_{1o}-1}_{j=m_{o}}2^{-jp'(s-(\frac{1}{p'}-\frac{1}{2})2\nu^*)}\II \left(  2^{j2\nu^*}>  \left[ \chi_{\alpha^*}(A)\right]^{\frac{2s_1}{2s+2\nu^*}-1} \right) \right)\nonumber\\
  &=& O \left( A^2 \left[ \chi_{\alpha^*}(A)\right]^{\frac{2s}{2s+2\nu^*}}\right).  \label{del'13}
\eeqn
\noindent Case 2: $s\leq(\frac{1}{p}-\frac{1}{2})2\nu^*$. Then, $d^*=\frac{2s'}{2s'+2\nu^*-1}$, and 
for $\Delta_{13}$ in \fr{del13}, one has 
\beqn  
  \Delta_{13}&=&O \left( A^2 \left[ \chi_{\alpha^*}(A)\right]^{1-\frac{p'}{2}}\sum^{J^*_{o}-1}_{j=m_{o}}2^{-j(p's'-(1-\frac{p'}{2})(2\nu^*-1))}\II \left(  2^{j2\nu^*}>  \left[ \chi_{\alpha^*}(A)\right]^{\frac{2s'}{2s'+2\nu^*}-1} \right)\right).\nonumber\\
  &=& O \left( A^2 \left[ \chi_{\alpha^*}(A)\right]^{1-\frac{p'}{2}}\sum^{J^*_{o}-1}_{j=m_{o}}2^{jp'((\frac{1}{p}-\frac{1}{2})2\nu^*-s)} \right)\nonumber\\
  &=& O \left( A^2 \left[ \chi_{\alpha^*}(A)\right]^{\frac{2s'}{2s'+2\nu^*-1}}\left[\ln(N)\right]^{\II\left(s=(\frac{1}{p}-\frac{1}{2})2\nu^*\right)}\right).  \label{del''13}
\eeqn
Note that the case $2^{J_1}< 2^{J^*_{o}}$ can be dealt with the same way, so we skip it.  \\
Finally, if $2^{J_1}< 2^{J_{o}}$, $\Delta_2$ can be decomposed into the following components
\beqn  
 \Delta_{21}&=&  O \left(\sum^{J-1}_{j=J_{0}}A^2 2^{-2js'} \right)= O\left(A^2 \left[\chi(A)\right]^d\right), \label{del21}\\
 \Delta_{22}&=&O \left( \sum^{J_{o}-1}_{j=J_{1}}\frac{ 2^{j_1(2\nu^*+1)}\ln^2(N)}{N}\II \left(  2^{j(2\nu^*+1)}\leq  \left[\chi(A)\right]^{d-1} \right)\right)=O\left(A^2 \left[\chi(A)\right]^d\right). \label{del22}\\
 \Delta_{23}&=&O \left( \sum^{J_{o}-1}_{j=J_{1}}A^{p'}2^{-p'js'}\left[\frac{\ln^2(N)}{N}2^{2\nu^* j}\right]^{1-p'/2}\II \left(  2^{j(2\nu^*+1)}>  \left[\chi(A)\right]^{d-1} \right)\right). \label{del23}
\eeqn
The derivation of the upper bounds for $\Delta_{23}$ in \fr{del23} will be dealt with the same way as \fr{del'13} and \fr{del''13}, so we skip it. Hence, combining \fr{e6}, \fr{ee14}, \fr{r21r31}, \fr{del11}, \fr{del12}, \fr{del'13}, \fr{del''13}, \fr{del21}, and \fr{del22}  completes the proof.  $\Box$

\end{document}